\documentclass[english]{article}
\usepackage{amssymb}
\usepackage{amsmath}
\usepackage[mathscr]{eucal}
\usepackage{mathrsfs}
 \usepackage{paralist}
\usepackage{color}

\newtheorem{theorem}{Theorem}[section]
\newtheorem{corollary}{Corollary}
\newtheorem{lemma}[theorem]{Lemma}
\newtheorem{proposition}{Proposition}

\newtheorem{definition}[theorem]{Definition}
\newtheorem{remark}{Remark}[section]

\newcommand{\Div}{\mbox{\rm div}\,}

\newcommand{\Int}[2]{{\displaystyle \int_{ #1}^{ #2}}}
\newcommand{\Lim}[1]{{\displaystyle \lim_{ #1}}}

\newcommand{\Max}[1]{{\displaystyle \max_{#1}}}

\newcommand{\Frac}[2]{\displaystyle{\frac{\displaystyle{#1}}{\displaystyle{#2}}}}

\newcommand{\beea}{\begin{eqnarray}}
\newcommand{\eeea}{\end{eqnarray}}

\newcommand{\bfe}{{\mbox{\boldmath $e$}} }

\newcommand{\bfq}{{\mbox{\boldmath $q$}} }
\newcommand{\bfz}{{\mbox{\boldmath $z$}} }
\newcommand{\0}{{\mbox{\boldmath $0$}} }

\newcommand{\BF}{\begin{footnotesize}}
\newcommand{\EF}{\end{footnotesize}}
\newcommand{\ode}[2]{{\displaystyle \frac{\mbox{$d #1$}}{\mbox{$d #2$}}}}

\newcommand{\bi}{\begin{itemize}}
\newcommand{\ei}{\end{itemize}}
\newcommand{\ed}{\end{document}}
\newcommand{\be}{\begin{equation}}
\newcommand{\ba}{\begin{array}}
\newcommand{\ea}{\end{array}}
\newcommand{\ee}{\end{equation}}
\newcommand{\eeq}[1]{\label{eq:#1}\end{equation}}
\newcommand{\real}{{\mathbb R}}

\newcommand{\nat}{{\mathbb N}}
\newcommand{\bfpsi}{\mbox{\boldmath $\psi$}}

\newcommand{\bfomega}{\mbox{\boldmath $\omega$}}

\newcommand{\bfx}{\mbox{\boldmath $x$}}

\newcommand{\bfv}{{\mbox{\boldmath $v$}} }
\newcommand{\bfu}{{\mbox{\boldmath $u$}} }
\newcommand{\bfw}{{\mbox{\boldmath $w$}} }

\newcommand{\bfa}{{\mbox{\boldmath $a$}} }

\newcommand{\bfp}{{\mbox{\boldmath $p$}} }

\newcommand{\bfA}{{\mbox{\boldmath $A$}} }

\newcommand{\bfR}{{\mbox{\boldmath $R$}} }

\newcommand{\bfG}{{\mbox{\boldmath $G$}} }

\newcommand{\bfH}{{\mbox{\boldmath $H$}} }
\newcommand{\bfI}{{\mbox{\boldmath $I$}} }

\newcommand{\bfh}{{\mbox{\boldmath $h$}} }

\newcommand{\cala}{{\cal A}}

\newcommand{\cald}{{\cal D}}

\newcommand{\calf}{{\cal F}}

\newcommand{\cali}{{\cal I}}

\newcommand{\cals}{{\cal S}}
\newcommand{\calt}{{\cal T}}

\newcommand{\bfT}{{\mbox{\boldmath $T$}} }
\newcommand{\bfV}{{\mbox{\boldmath $V$}} }

\newcommand{\bfW}{{\mbox{\boldmath $W$}} }
\newcommand{\bfF}{{\mbox{\boldmath $F$}} }
\newcommand{\bfK}{{\mbox{\boldmath $K$}} }
\newcommand{\bfL}{{\mbox{\boldmath $L$}} }

\newcommand{\bfg}{{\mbox{\boldmath $g$}} }
\newcommand{\bfn}{{\mbox{\boldmath $n$}} }
\newcommand{\bfs}{{\mbox{\boldmath $s$}} }

\newcommand{\half}{\mbox{$\frac{1}{2}$}}

\def\Bbb R{\real}

\def\tilde{\widetilde}
\def\bar{\overline}

\newcommand{\ED}{\end{description}}

\newcommand{\dist}{\mbox{\rm dist\,}}

\newcommand{\Br}{\begin{remark}\begin{rm}}
\newcommand{\Er}{\end{rm}\end{remark}}
\newcommand{\Bt}{\begin{theorem}\begin{sl}}
\newcommand{\Et}{\end{sl}\end{theorem}}
\newcommand{\Bl}{\begin{lemma}\begin{sl}}
\newcommand{\El}{\end{sl}\end{lemma}}
\newcommand{\Eqref}[1]{{\rm (\ref{eq:#1})}}
\newcommand{\Bc}{\begin{corollary}\begin{sl}}
\newcommand{\Ec}{\end{sl}\end{corollary}}
\newcommand{\ET}[1]{\end{sl}\label{theorem:#1}\end{theorem}}
\newcommand{\EL}[1]{\end{sl}\label{lemma:#1}\end{lemma}}
\newcommand{\theoref}[1]{{\rm Theorem \ref{theorem:#1}}}
\newcommand{\ER}[1]{\end{rm}\label{remark:#1}\end{remark}}
\newcommand{\EC}[1]{\end{sl}\label{corollary:#1}\end{corollary}}
\newcommand{\remref}[1]{{\rm Remark \ref{theorem:#1}}}
\newcommand{\cororef}[1]{{\rm Corollary \ref{corollary:#1}}}
\newcommand{\lemmref}[1]{{\rm Lemma \ref{lemma:#1}}}
\newcommand{\essup}[1]{{\rm ess}\,{{\displaystyle \sup_{\hspace*{-6mm}{#1}}}}\,}
\renewcommand{\nat}{{\mathbb N}}

\newcommand{\propref}[1]{{\rm Proposition \ref{proposition:#1}}}
\newcommand{\Bp}{\begin{proposition}\begin{sl}}
\newcommand{\EP}[1]{\end{sl}\label{proposition:#1}\end{proposition}}

\newcommand{\QED}{\par\hfill$\square$\par}
\usepackage{graphicx}
\graphicspath{%
    {converted_graphics/}
    {/}
}
\begin{document}
\pagenumbering{arabic}
\title{Inertial Motions of a Rigid Body with a cavity filled with a viscous liquid} 
\medskip\bigskip
\author{Giovanni P. Galdi,\ Giusy Mazzone \& Paolo Zunino}
\date{}
\maketitle
\newcommand{\vp}{\varphi}
\renewcommand{\theequation}{\arabic{section}.\arabic{equation}}
\begin{abstract}
We study inertial motions of the coupled system, $\mathscr S$, constituted by a rigid body containing a cavity that is completely filled with a viscous liquid. We show that for data of arbitrary size (initial kinetic energy and total angular momentum) every weak solution (\`a la Leray--Hopf) converges, as time goes to infinity, to a uniform rotation, thus proving a famous ``conjecture" of {\sc N.Ye. Zhukovskii}. Moreover we show that, in an appropriate range of initial data, this rotation must occur along the central axis of inertia of $\mathscr S$ with the larger moment of inertia. Furthermore, we provide necessary and sufficient conditions for the rigorous nonlinear stability of permanent rotations, which improve and/or generalize results previously given by other authors under different types of approximation of the original equations, and/or suitable symmetry assumptions on the shape of the cavity. Finally, we present a number of results obtained by a targeted numerical simulation that, on the one hand, complement the analytical findings, whereas, on the other hand, point out new features that the analysis is yet not able to catch, and, as such, lay the foundation for interesting and challenging future investigation. 
\end{abstract}
\setcounter{section}{0}
\section*{1. Introduction~\footnote{Some of the main results proved in this paper were previously announced in \cite{GMZ}.}}
Consider a  body $\mathscr B$ (rigid or deformable) freely moving in the physical space. We assume that no {\em external} forces are acting on $\mathscr B$, so that $\mathscr B$ executes an {\em inertial motion}. 

In an inertial motion, the center of mass ${G}$ of $\mathscr B$  moves by uniform and rectilinear motion. Consequently, the relevant dynamics of $\mathscr B$  reduces to the study of its movement about ${G}$, namely, with respect to a  frame having  origin in ${G}$ and axes of constant direction in a given inertial frame, $\cali$. \renewcommand{\theequation}{1.\arabic{equation}}\par
\smallskip\par
Assume first  $\mathscr B$  {\em rigid}. Then, the dynamics about ${G}$ is governed  by the law of conservation of the total angular momentum with respect to ${G}$, $\bfK_{{G}}$.  As is well known,  this law is equivalent to the  classical Euler equation, written in  the body-fixed central frame of inertia with origin in ${G}$ and  axes ({\em central axes of inertia}) parallel to three orthogonal eigenvectors $\bfe_i$, $i=1,2,3,$ of the inertia tensor $\bfI$ of $\mathscr B$ evaluated with respect to ${G}$ ({\em central inertia tensor}), namely,
\be
\bfI\cdot\dot{\bfomega}+\bfomega\times\bfI\cdot\bfomega=\0\,, 
\eeq{01}
with $\bfomega$  angular velocity of $\mathscr B$.
The dynamics of $\mathscr B$ about ${G}$, obtained by solving \Eqref{01},  is very rich.  In particular, time independent  motions ({\em permanent rotations}) may occur if and only if the (constant) angular velocity $\bfomega$  is directed along one of the $\bfe_i$, which, in turn, must align with $\bfK_{{G}}$. On the opposite side, time-dependent  motions may be very complicated depending on the ``symmetry'' of the distribution of mass of $\mathscr B.$ In mathematical terms, this  is expressed by the properties of the eigenvalues, $A$, $B$, and $C$, of $\bfI$ corresponding to $\bfe_1$, $\bfe_2$, and $\bfe_3$, respectively. For example, if $A= B\neq C$, then the most general  motion of $\mathscr B$  about $G$ is a {\em regular precession}, where $\mathscr B$ rotates uniformly around the axis parallel to $\bfe_3$ and passing through $G$, while the latter rotates uniformly around the direction of $\bfK_G$. However, if $A\neq B\neq C$, the generic motion can be quite involved, and falls in the class of the {\em motions \`a la Poinsot}; see, e.g., \cite[Chapter 1, \S4]{Lei}. 

Assume, next, that inside  $\mathscr B$ we perform a hollow cavity $\mathscr C$, completely contained in $\mathscr B$, and  fill it up with a {\em viscous} liquid $\mathscr L$. Then,  the inertial motions of the {\em coupled} system $\mathscr S\equiv\mathscr B\cup\mathscr L$ are expected to show, at large times, an altogether  peculiar behavior. In fact, after an initial  time interval --whose length depends on the size of the initial data and the viscosity of the liquid--  where the motion is typically of chaotic nature (especially for ``small'' viscosity, see \cite{LK} and Section 9.1), the coupled system reaches eventually a more orderly configuration,  due to the combined effects of viscosity and incompressibility  of the liquid \cite{Zu,MoRu,Ch}. 
This circumstance was first pointed out, back in 1885, by N.Ye. Zhukovskii who, under the assumption that $\mathscr L$ is a Navier--Stokes liquid,  in \cite[p. 152]{Zu} stated the following property: ``{\em  The motions of $\mathscr S$ about ${G}$ will eventually 
be rigid motions and, precisely, permanent rotations,  no matter the size and shape of $\mathscr C$, the viscosity of $\mathscr L$, and the initial movement of $\mathscr S$}''. 

Such an intriguing property  is quite in contrast with the case, described earlier on, when the cavity is empty, and can be viewed as a {\em stabilizing effect} of the liquid on the motion of the solid.  Zhukhovskii's argument is simply 
described as follows. Because of viscous effects, the velocity field of $\mathscr L$ relative to $\mathscr B$ must eventually vanish, so that $\mathscr S$ will eventually move by rigid motion. Under this condition, from the equation of linear momentum of the liquid we derive that the pressure gradient in $\mathscr L$, $\nabla {\sf p}$, must balance centrifugal forces:
\be
\rho\,[\dot{\bfomega}\times\bfx+\bfomega\times(\bfomega\times\bfx)]=\nabla{\sf p}\,,
\eeq{zaza}  
where 
$\rho$ is the density of $\mathscr L$. Therefore, by taking the curl of both sides of the equation we get $\dot{\bfomega}={\bf 0}$, which shows that only permanent rotations may occur. 

Clearly, the previous reasoning is purely heuristic and takes for granted the following two facts: (i) the liquid comes to rest (relative to a body-fixed frame) as time goes to infinity, and (ii) once (i) occurs, the dynamics of the system is governed by \Eqref{zaza}. However, neither of these properties is obvious from a rigorous  viewpoint. For one thing, even   
though the {\em total} energy is decreasing, it is not said that the velocity of the liquid relative to the body must eventually vanish, especially when the initial data are such that the magnitude of the total angular momentum (that is {\em conserved} throughout the motion) is arbitrary large. 
Furthermore, admitting property (ii) would mean, in terms of dynamical systems, that the  
$\Omega$-limit set\footnote{Throughout this paper, we use the upper case notation $(\Omega$) in place of the usual lower case one ($\omega$) to avoid confusion with the angular velocity.} of a generic trajectory is {\em not empty} and, more importantly, {\em invariant} in the relevant class of solutions, a statement that, also when properly formulated (see further on), is far from being evident.

It must be emphasized that, even though never rigorously proved right or wrong, Zhukovskii's statement was later on elevated to the rank of a theorem, mainly by Russian authors; see, e.g., \cite[p. 98]{MoRu}, \cite[p. 3]{Ch}. We prefer instead to refer to it  as {\em Zhukhovskii's conjecture}.

In this respect, we would like to observe that {there is an} abundant mathematical literature dedicated to the motion of a solid with an interior cavity filled with a viscous liquid, that ranges from the early contributions of Stokes \cite{Stokes0}, Zhoukhovski \cite{Zu}, and Hough \cite{Hough}  to the more recent classical  monographs \cite{MoRu,Ch,KK}. {However}, results are rarely of an exact nature, either due to the approximate models or else due to an approximate mathematical treatment.

The primary goal of this paper is to  perform a comprehensive and rigorous mathematical analysis of  the  inertial motions of a rigid body with an interior cavity completely filled by a viscous liquid, with special regard to  their asymptotic behavior in time. 

Our objective is, in fact, manifold and is presented next. 

In the first place, we investigate the validity of Zhukovskii's conjecture. This study is performed  in a very general class, $\mathcal S$, of motions  described by appropriate weak solutions \`a la Leray-Hopf.  Such solutions  are characterized by the velocity field, $\bfv$, of the liquid $\mathscr L$ relative to $\mathscr B$, and angular velocity, $\bfomega$, of $\mathscr B$  satisfying (in a suitable sense)  balance of total energy and conservation of total angular momentum at all times. The class $\cals$ is   not empty if only the initial data have a finite (total) energy.  We then prove that   each solution $(\bfv,\bfomega)$ in $\cals$, as $t\to\infty$, must converge (in proper topology) to a state of motion of $\mathscr S$ where $\bfv\equiv \0$ and $\bfomega=\bar\bfomega$, for some constant vector $\bar \bfomega$. This means that, eventually, $\mathscr S$ moves as a single rigid body by constant rotational motion, which is exactly Zhukhovskii's conjecture. 

In the above argument,  the proof of convergence represents, as expected, the most challenging issue, and we  carry it out along the following steps. For a given weak solution, we begin to show that the corresponding $\Omega$-limit set is not empty and of the form $\{\bfv\equiv \0\}\times \cala$, where $\cala$ is a compact, connected subset of $\real^3$ that is left invariant by the semigroup associated to the Euler equations \Eqref{01}, with $\bfI$ central inertia tensor of $\mathscr S$; see \propref{3}. From the physical viewpoint, this means that, ultimately, the system  $\mathscr S$ moves, as a whole, by rigid motion, but, at this stage, it is not said that such a motion is a permanent rotation. However, it is  readily seen that the latter certainly holds if the $\Omega$-limit set is invariant in the class $\cals$ of {\em weak solutions}; see \lemmref{5}.  Thus, the proof of the conjecture boils down to the proof of the stated invariance of the $\Omega$-limit set. Now, as is well known,   the invariance property requires, typically, the {\em uniqueness} of solutions (see, e.g., \cite{DaSl,Hale}), which unfortunately is not available in the case at hand, in that our weak solutions possess (in their liquid variable) the same  properties as the three-dimensional weak solutions to the Navier--Stokes equations. Nevertheless, due to the circumstance that the velocity field of the liquid tends to zero asymptotically, we  prove that, in fact,  uniqueness for ``large'' times only would suffice. Therefore, the question arises of whether any weak solution becomes ``strong'' (and therefore unique) for sufficiently large times. In the case of the Navier--Stokes equations this is true (and well established), since  every weak solution\footnote{Corresponding to vanishing external force and boundary data.}  becomes ``small'' as time goes to infinity, regardless of the ``size'' of the initial data. The same property is by no means obvious in the case at hand, due to the
presence of an, in general, large ``conservative'' component (conservation
of the total angular momentum).
One may then guess the property to be true for ``small'' total angular
momentum, but what for arbitrarily large? This non-trivial question is addressed and solved in \propref{2}, where it is shown that, if the cavity is sufficiently smooth (e. g., of class $C^2$) for {\em any} weak solution there exists a corresponding time, after which the solution becomes strong and, in particular, unique. With the help of this result, we then prove the invariance of the $\Omega$-limit set; see \propref{4}. Once all the above properties have been secured, we are then able to show our first main achievement in \theoref{3}. In this theorem,  we provide a rigorous proof of Zhukovskii's conjecture, as well as we  show that the permanent rotation must occur along one of the central axes of inertia of $\mathscr S$, and that the latter must eventually be aligned with the direction of the (constant) total angular momentum of $\mathscr S$ with respect to $G$, $\bfK_G$.

It must be noticed that \theoref{3} does not specify around which axis the ``final'' permanent rotation will be attained. This issue is of great significance also because  our weak solutions may lack of uniqueness, and consequently we may have, in principle, two different solutions with the same initial data generating, asymptotically, two permanent rotations around different axes.  The other main objective of this paper is therefore to analyze this problem in some details. In particular, we shall prove that, for initial data in an  appropriate range, the   permanent rotation will always occur along that central axis, ${\sf a}$, with the larger moment of inertia; see \theoref{4} and \remref{4}. This phenomenon may have a natural physical explanation. In fact, the coupled system $\mathscr S$, moving in absence of external forces,  tends eventually to reach the state of ``minimal motion". Because of conservation of total angular momentum, the latter  cannot be the rest, unless the initial total angular momentum is zero, a case that is of no intrinsic relevance; see Remark \ref{rem:ST}.  Consequently $\mathscr S$ chooses, in general, to rotate around the axis where the spin is a (non-zero) minimum.
In addition to the property just described, we also show that, provided the two components of the angular velocity in  direction orthogonal to {\sf a} are {\em initially} sufficiently ``small'', and the {\em initial} relative velocity of the liquid is zero, the permanent rotation is attained with ${\sf a}$ and $\bfK_G$  keeping the {\em same} orientation they had at time $t=0$. Namely, if the angle, $\theta(t)$, between ${\sf a}$ and $\bfK_G$ is such that $\theta(0)<\pi/2$, then $\theta(\infty)=0$, whereas if $\theta(0)>\pi/2$, then $\theta(\infty)=\pi$; see Remark  \ref{rem:6}. Notice that the same property may fail to hold if the ``smallness'' assumption is violated, as suggested by the numerical test presented in Section 9.2; see Figure 5, top panel.   

A related, and final, goal of our analytic investigation concerns the stability of permanent rotations. In this connection, 
we recall that the study of the stability of  permanent rotations of a body with an interior cavity completely filled with a liquid has a long history, beginning with  the work of Sobolev \cite{Sob} (in the inviscid case),  Rumyantsev \cite[p. 203]{Rum} (see also \cite[\S 3-3]{MoRu}), Chernousko \cite{Ch}, and Smirnova \cite{Smirn1,Smirn2} until the more recent work of Kostyuchenko {\em et al.} \cite{KSY}, and Kopachevsky \& Krein \cite{KK}.\footnote{In the latter two cited references one can find also a detailed account of previous relevant contributions to the stability problem.}  With the exception of \cite{Rum,MoRu}, 
results proved by these authors are obtained under different simplifications. More precisely, in \cite{KSY,KK}  stability/instability is investigated by linearizing the relevant equations. However, the corresponding results need not be valid for the original nonlinear problem due to the lack, to date, of a ``linearization principle'' that may validate the above findings  at the nonlinear level.  Similar results are obtained  in \cite{Ch,Smirn1,Smirn2} under the assumption of large/small viscosity (small/large Reynolds number) and/or geometrical symmetry of the cavity, which allows one to approximate  the original problem with  an appropriate system of nonlinear {\em ordinary} differential equations.  Finally, in  \cite{Rum} a nonlinear stability analysis is performed that provides sufficient conditions for stability of permanent rotations.\footnote{As a matter of fact, the analysis in \cite{Rum} is still {\em formal}, because it lacks of a suitable corresponding existence theorem for the relevant equations.} In contrast to the above results, our approach to the stability problem is rigorous, on the one hand, and very straightforward on the other hand. Actually, it is based upon the attainability results of \theoref{4} along with balance of total energy and conservation of total angular momentum. In this way we are able  to show {\em necessary and sufficient} conditions for stability (in the sense of Lyapunov) for the full nonlinear problem,  {\em without any approximation or assumptions on the shape of the cavity}. These conditions contain those of  \cite{Rum} as a particular case, and extend those of \cite{KSY,Ch,Smirn1,Smirn2} to the nonlinear level; see \theoref{5}.

The analytic study just described is followed and enriched by a number of results obtained from suitably targeted  numerical tests; see Section 9.  These tests have a multiple purpose. In the first place, they suggest a {\em quantitative} estimate of the critical time, $t_c$, after which the coupled system executes, within a given tolerance, the ultimate  permanent rotation. All other data being fixed, it is found that  $t_c\sim (\nu)^{-0.3}$, where $\nu$ is the coefficient of kinematic viscosity of the liquid; see Section 9.1. In the second place, they complement  the analytical findings obtained in \theoref{4}, by investigating   whether the restrictions on the range of the initial data are indeed sharp, or can be possibly removed and/or improved; see Section 9.2. Finally, they point out another interesting phenomenon that the mathematical analysis is yet not able to catch. As we have described earlier on, the final permanent rotation occurs around  the central axis {\sf a} with the largest moment of inertia, and, at least under suitable ``smallness'' assumptions,  having same or opposite orientation as $\bfK_G$, depending on whether the {\em initial} angle $\theta(0)$ between ${\sf a}$ and $\bfK_G$ is less or greater than $\pi/2$. Now, if $\theta(0)=\pi/2$,  numerical tests show a ``flip-over'' phenomenon triggered by the viscosity of the liquid. Precisely, as the kinematic viscosity coefficient $\nu$ is decreased from a ``sufficiently large'' value --all other data being kept fixed-- we find a critical value $\nu_c$ such that the orientation of ${\sf a}$ and $\bfK_G$ are the same or opposite according to whether $\nu>\nu_c$ or $\nu<\nu_c$. The analytical investigation of this phenomenon will be  object  of  future work.
\renewcommand{\theequation}{2.\arabic{equation}}
\setcounter{section}{1}
\section{Notation and Preparatory Results}
The notation we shall use throughout is  quite standard. The symbol $\nat$ denotes the set of positive integers. By $\real$ we denote the set of real numbers, and by $\real^3$ the Euclidean three-dimensional space. Vectors in $\real^3$ will be indicated by boldfaced letters, and the  canonical base, $B$, in $\real^3$ by $\{\bfe_1,\bfe_2,\bfe_3\}$. Components of a vector $\bfv$ in $B$ are indicated by $(v_1,v_2,v_3)$, whereas $|\bfv|$ represents the magnitude of $\bfv$. \par  
If $O$ is a point in $\real^3$ and $\bfe$ is a unit vector, by $\{O,\bfe\}$ we mean the straight line parallel to $\bfe$ and passing through $O$.
\par
Let $A$ be an open set of $\real^3$. We denote by $L^2(A)$, and $W^{k,2}(A), W_0^{k,2}(A)$, $k\in \nat$, the usual Lebesgue and Sobolev spaces, with norms $\|\cdot\|_2$ and $\|\cdot\|_{k,2}$, respectively.\footnote{Unless confusion arises, we shall use the same symbol for spaces of scalar, vector and tensor functions. Moreover, in  integrals we shall typically omit the infinitesimal element of integration.} 

If $\{\bfG,\bfH\}$, $\{\bfg,\bfh\}$, $\{g,h\}$ are pairs of second-order tensor,  vector and scalar fields on $A$, respectively, we set $$
(\bfG,\bfH)_A=\int_A G_{ij}H_{ij}\,, \ \ \
(\bfg,\bfh)_A=\int_A g_ih_i\,,\ \ \ (g,h)_A=\int_A g\,h\,,
$$
where, typically, we shall omit the subscript $``A"$.

If $A$ is a bounded, Lipschitz 
domain with outward unit normal $\bfn,$ it is well known that the Helmholtz-Weyl decomposition holds (e.g., \cite{Gabook}): $L^2 (A)= H(A)\oplus G(A),$
 with
\[
H(A) = \{\bfu \in L^2 (A) : \Div \bfu =0 \textrm{ and } \bfu \cdot \bfn{|_{\partial A}} = 0 \}\,,
\]
where $\Div\bfu$ and $\bfu \cdot\bfn|_{\partial A}$ are understood in the sense of distributions, and
\[
G(A)= \{\bfw \in L^2 (A) : \bfw = \nabla p, \textrm{ for some } p \in W^{1,2}( A)\}.
\]
The projection operator from $L^2(A)$ onto $H(A)$ is denoted by $\mathscr P.$ 

We also set $\cald_0^{1,2}(A):=H(A)\cap W_0^{1,2}(A)$.

If $X$ is a Banach space with norm $\|\cdot\|_X$, and $I\subset\real$ an interval, we denote by $L^q(I;X)$ [resp. $W^{k,q}(I;X)$, $k\in\nat$], the space of functions $f:I\mapsto X$ such that
$  
\left(\int_I\|f(t)\|_X^qdt\right)^{1/q}<\infty
$
[resp. $\sum_{l=0}^k\left(\int_I\|\partial_t^{l}f(t)\|_X^qdt\right)^{1/q}<\infty$]. Likewise, we write $f\in C^k(I;X)$ if $f$ is $k$-times differentiable with values in $X$ and $\max_{t\in I}\|\partial_t^lf(t)\|_X<\infty$, $l=0,1,\ldots,k$.
Moreover, $f\in C_w(I;X)$ means that the map $t\in I\to \ell(f(t))\in\real$ is continuous for all bounded linear functionals $\ell$ defined on $X$.  
 
In denoting all the above spaces, we shall omit the symbol $X$ if $X=\real$. 
\par
We conclude by noting that with the symbols $c$, $c_0$, $c_1$, etc.,  we denote positive constants, whose particular value is unessential to the context. 
\smallskip\par
We shall now collect some results of Gronwall's lemma-type that will be referred to in the paper. Even though they can be probably found in the existing literature, nevertheless we shall sketch their proof in the Appendix  for completeness. 
\Bl 
Suppose $y\in L^\infty(0,\infty)$, $y\ge 0$, satisfies for a.a. $s\ge 0$ and all $t\ge s$,
\be 
y(t)\le y(s)-k\int_s^ty(\tau)\,d\tau+\int_s^tF(\tau)\,d\tau\,,
\eeq{2.1}
where $k>0$, whereas $F\in L^q(a,\infty)\cap L^1_{{\rm loc}}(0,\infty)$ for some $a>0$,  $q\in [1,\infty)$, and $F(t)\ge 0$, a.a. $t\ge 0$. Then
$$
\lim_{t\to\infty}y(t)=0\,.
$$ 
Moreover, if in particular $F\equiv 0$, then
$$
y(t)\le y(s)\,{\rm e}^{-k(t-s)}
\,,\ \ \mbox{for all $t\ge s$}\,.
$$
\EL{1}

\Bl Let $y:[t_0,t_1)\to[0,\infty)$,  $t_1>t_0\ge 0$, be an absolutely continuous function satisfying for some $a,b,c,\delta>0$\footnote{For the sake of completeness, we would like to remark that, as is well known, the lemma continues to hold if $c=0$, even without assuming the first condition in (ii).} and $\alpha>1$,
\begin{itemize}
\item[{\rm (i)}]\ $y'\le -a\,y+b\,y^\alpha+c$\ \ \mbox{in $(t_0,t_1)$}\,;
\item[{\rm (ii)}]\ $\Int{t_0}{t_1}y(\tau)\,d\tau<\Frac{\delta^2}{4c}\,,\ \ y(t_0)<\Frac{\delta}{\sqrt{2}}$\,. 
\end{itemize}
Then, if $k:=-a+b\,\delta^{\alpha-1}<0$, we have
\be
y(t)<\delta\,,\ \ \mbox{for all $t\in [t_0,t_1)$}\,.
\eeq{2.7}
Moreover, if $t_1=\infty$ we have also
\be
\lim_{t\to\infty}y(t)=0\,.
\eeq{2.8}
\EL{2}\renewcommand{\theequation}{3.\arabic{equation}}
\setcounter{equation}{0}
\section{Formulation of the Problem and Preliminary Considerations} Let $\mathscr B$ denote a rigid body with a cavity, $\mathscr C$, completely contained in it. In mathematical terms, $\mathscr B:=\Omega_1\backslash\bar{\Omega_2}$, $\mathscr C:=\Omega_2$, where $\Omega_i$, $i=1,2$ are simply connected, bounded domains of $\real^3$ with $\bar{\Omega_2}\subset\Omega_1$. Throughout the paper, {\em we shall assume that $\mathscr C$ is of class} $C^2$, even though some peripheral results continue to  hold under weaker (or even without any)  hypothesis. 

We suppose that $\mathscr C$ is completely filled with a viscous liquid whose motion is governed by the Navier--Stokes equations. We also assume that the coupled system $\mathscr S:=\mathscr B\cup\mathscr C$ is {\em isolated}, namely, no external forces are applied, so that $\mathscr S$ moves by {\em inertial motion}.  

Under these premises, the center of mass, $G$, of $\mathscr S$ will move by uniform and rectilinear motion with respect to an inertial frame, $\cali$, so that the relevant dynamics of $\mathscr S$ reduces to the study of its motion {\em about $G$}, namely, with respect to an (inertial) frame  with the origin at $G$ and axes whose direction is constant in $\mathcal I$. However, following the approach of classical mechanics of rigid bodies, the motion of $\mathscr S$ is better described in the (non-inertial) frame, $\calf$, attached to $\mathscr B$, and having origin at $G$ and axes directed along the principal axes, $\{\bfe_1,\bfe_2,\bfe_3\}$,  of the central inertia tensor\footnote{Namely, the inertial tensor with respect to $G$.} $\bfI$ of $\mathscr S$.  Notice that with respect to $\mathcal F$, the domain, $\mathscr C$, occupied by the liquid is independent of time. Thus,  when referred to $\mathcal F$, the relevant equations of motion of  $\mathscr S$ are given by (see \cite{MoRu}, Chapter 1; see also \cite{Mazz}, Chapter 1) 
\be\ba{ll}\smallskip\left.\ba{cc}\smallskip
\rho(\partial_t\bfu+\bfv\cdot\nabla\bfu+\bfomega
\times\bfu)=\mu\Delta\bfu-\nabla {\sf p}\\
\Div\bfu=0\ea\right\}\ \ \mbox{in $\mathscr C\times (0,\infty)$}\\
\bfI_{\mathscr B}\cdot\ode{\bfomega}{t}+\bfomega\times (\bfI_{\mathscr B}\cdot\bfomega)=-\Int{\partial\mathscr C}{}\bfx\times \bfT(\bfu,{\sf p})\cdot\bfn\,,\ \ \mbox{in $(0,\infty)$}\,.
\ea
\eeq{3.1}
Here, $\bfu, {\sf p}$ are absolute velocity and pressure fields of the liquid, $\rho$ is its (constant) density,  $\mu$ its shear viscosity coefficient, and $\bfn$ is the unit outer normal to $\partial\mathscr C$. Moreover, $\bfomega$ and $\bfI_{\mathscr B}$ are angular velocity and inertia tensor of $\mathscr B$ with respect to $G$, respectively. Finally, $\bfv:=\bfu-\bfomega\times\bfx$ is the relative velocity field of the liquid, and $\bfT=\bfT(\bfw,{\sf \phi})$ is the Cauchy stress tensor defined as
\be
\bfT(\bfw,{\sf \phi})=\mu(\nabla\bfw+(\nabla\bfw)^\top)-{\sf \phi}\bf1,
\eeq{3.2}
with $\bf1$ unitary tensor. To \Eqref{3.1} we append the no-slip condition of the liquid at the bounding wall of the cavity
\be
\bfu=\bfomega\times\bfx\,,\ \ \mbox{at $\partial\mathscr C\times (0,\infty)$}\,.
\eeq{3.3}\par
Let 
\be
\bfA:= \bfI\cdot\bfomega+\Int{\mathscr C}{}\rho\,\bfx\times \bfv
\eeq{3.4}
denote the {\em total angular momentum} of  $\mathscr S$ with respect to $G$, referred to $\calf$.
In the Appendix it is then shown that the system of equation \Eqref{3.1}--\Eqref{3.3} is {\em equivalent} to the following one
\be\ba{ll}\smallskip\left.\ba{cc}\smallskip
\rho(\partial_t\bfv+\bfv\cdot\nabla\bfv+\dot{\bfomega}\times\bfx+2\bfomega
\times\bfv)=\mu\Delta\bfv-\nabla {\rm p}\\
\Div\bfv=0\ea\right\}\ \ \mbox{in $\mathscr C\times (0,\infty)$}\\
\ode{\bfA}{t}=\bfA\times\bfomega\,,\ \ \mbox{in $(0,\infty)$}\,,
\ea
\eeq{3.5}
along with the boundary condition
\be
\bfv=\0\,,\ \ \mbox{at $\partial\mathscr C\times (0,\infty)$}\,,
\eeq{3.6}
where in  \Eqref{3.5} the dot indicates differentiation with respect to time, and ${\rm p}={\sf p}+\half (\bfomega\times\bfx)^2$. 

The form \Eqref{3.5} of the relevant equations provides a better insight of the mixed character of the problem. Actually,  equations \Eqref{3.5}$_{1,2}$ describe the ``dissipative component'' of $\mathscr S$, while  \Eqref{3.5}$_{3}$ characterizes its ``conservative component''. As a matter of fact, \Eqref{3.5}$_{3}$ is just the equation of  conservation of total angular momentum written in the non-inertial frame $\calf$. 

Our expectation is that, as time goes to infinity, the coupled system $\mathscr S$ will move, as a whole, by  rigid motion. Of course, the angular velocity of this ``final'' motion must be compatible with conservation of the total angular momentum of $\mathscr S$. With this in mind, we introduce the (unknown) ``final'' angular velocity
\be   
\bfomega_\infty:=\bfI^{-1}\cdot\bfA\equiv \bfomega+\bfI^{-1}\cdot\left(\Int{\mathscr C}{}\rho\,\bfx\times \bfv\right)
\eeq{3.7}
and rewrite \Eqref{3.5} equivalently in the following form
\be\ba{cc}\medskip\left.\ba{rl}\smallskip
\rho(\partial_t\bfv+\bfv\cdot\nabla\bfv+(\dot{\bfomega}_\infty+\dot{\bfa})\times\bfx+2(\bfomega_\infty&\!+\bfa)
\times\bfv)\\ \medskip&=\mu\Delta\bfv-\nabla {\rm p}\\
\Div\bfv&=0\ea\right\}\ \ \mbox{in $\mathscr C\times (0,\infty)$}\\
\bfI\cdot\dot{\bfomega}_\infty+(\bfomega_\infty+\bfa)\times(\bfI\cdot\bfomega_\infty)=\0\,,\ \ \mbox{in $(0,\infty)$}\,,
\ea
\eeq{3.8}
where
\be
\bfa:=-\rho\bfI^{-1}\cdot\left(\int_{\mathscr C}\bfx\times\bfv\right)\,.
\eeq{3.9}
\par
The main results of this paper are primarily  based on the study of the asymptotic behavior in time of pairs $(\bfv,\bfomega_\infty)$ satisfying \Eqref{3.8}--\Eqref{3.9} for a suitable ${\rm p}$,  and  side condition \Eqref{3.6}. This investigation will be carried out in a very broad class of weak solutions (\`a la Leray-Hopf) that will be introduced in the next section. An important tool to this effect is provided by the equation of balance of the total energy of $\mathscr S$. Before deriving (formally) the latter, we premise the following result due to \cite[\S\S 7.2.2--7.2.4]{KK} that will turn out to be useful  in several circumstances.    
\Bl Let $\bfw\in H(\mathscr C)$. 
Then, there is $c>0$ such that
$$
\|\bfw\|_2^2\ge \|\bfw\|_2^2-\left(\rho\,\bfI^{-1}\cdot\int_{\mathscr C}\bfx\times\bfw\right)\cdot\left(\int_{\mathscr C}\bfx\times\bfw\right)\ge c\|\bfw\|_2^2\,.
$$
\EL{3}

In order to derive the balance of energy equation, we dot-multiply both sides of \Eqref{3.8}$_1$ by $\bfv$ and integrate by parts over $\mathscr C$ to obtain, with the help of \Eqref{3.8}$_2$ and \Eqref{3.6},  
\be
\frac{\rho}{2}\,\ode{}{t}\|\bfv\|_2^2+\mu\|\nabla\bfv\|_2^2=-\dot{\bfomega}_\infty\cdot\int_{\mathscr C}\rho\bfx\times\bfv-\dot{\bfa}\cdot\int_{\mathscr C}\rho \bfx\times\bfv\,.
\eeq{3.10}
Now, recalling that $\bfI$ is a symmetric tensor, from \Eqref{3.8}$_3$ and \Eqref{3.9} we infer
\be
-\dot{\bfomega}_\infty\cdot\int_{\mathscr C}\rho\,\bfx\times\bfv=\bfa\cdot\bfI\cdot\dot{\bfomega}_\infty=-\bfomega_\infty\times
(\bfI\cdot\bfomega_\infty)\cdot\bfa\,.
\eeq{3.11}
Furthermore, by dot-multiplying both sides of \Eqref{3.8}$_3$ by $\bfomega_\infty$ we get
$$
\frac{1}{2}\ode{}{t}(\bfomega_\infty\cdot\bfI\cdot\bfomega_\infty)=\bfomega_\infty\times(\bfI\cdot\bfomega_\infty)\cdot\bfa\,,
$$
whereas, by \Eqref{3.9},
$$
-\dot{\bfa}\cdot\int_{\mathscr C}\rho \bfx\times\bfv=\frac{1}{2}\ode{}{t}
(\bfa\cdot\bfI\cdot\bfa)\,.
$$
Combining the last two displayed equations with \Eqref{3.10} and \Eqref{3.11} we conclude
\be
\ode{\mathcal E}{t}=-\mu\|\nabla\bfv\|_2^2
\eeq{3.12}
where
\be
\mathcal E:=\frac{1}{2}\left(\rho\,\|\bfv\|_2^2-\bfa\cdot\bfI\cdot\bfa+\bfomega_\infty\cdot\bfI\cdot\bfomega_\infty\right)\,.
\eeq{3.13}
Equation \Eqref{3.12} represents the balance of energy for the coupled system $\mathscr S$.
Notice that $\mathcal E$ is positive definite in the pair $(\bfv,\bfomega_\infty)$   since, by \lemmref{3},
\be
c_1\|\bfv\|_2^2\le E:=\frac{1}{2}\left(\rho\,\|\bfv\|_2^2-\bfa\cdot\bfI\cdot\bfa\right)\le c_2\|\bfv\|_2^2\,.
\eeq{3.14}\renewcommand{\theequation}{4.\arabic{equation}}
\setcounter{equation}{0}
\section{Weak Solutions and their Basic Properties}
We now  introduce a suitable class of weak solutions to the problem \Eqref{3.6}, \Eqref{3.8}--\Eqref{3.9}. To this end, dot-multiplying both sides of \Eqref{3.8}$_1$ by $\bfpsi\in\cald_0^{1,2}(\mathscr C)$, and  integrating by parts over $\mathscr C\times (0,t)$, we deduce
\be\ba{rl}\smallskip
(\rho\,\bfv(t),\bfpsi)+&\rho(\bfomega_\infty(t)+\bfa(t))\cdot\Int{\mathscr C}{}\bfx\times\bfpsi\\ \smallskip
&+\Int{0}t\left\{\rho(\bfv\cdot\nabla\bfv,\bfpsi)+2\rho((\bfomega_\infty+\bfa)\times\bfv,\bfpsi)+\mu(\nabla\bfv,\nabla\bfpsi)\right\}\\ \smallskip
=&(\rho\,\bfv(0),\bfpsi)+\rho(\bfomega_\infty(0)+\bfa(0))\cdot\Int{\mathscr C}{}\bfx\times\bfpsi\,,\\
&\qquad\qquad\qquad\mbox{for all $\bfpsi\in\cald_0^{1,2}(\mathscr C)$ and all $t\in (0,\infty)$.}
\ea
\eeq{4.1}
Moreover, integrating \Eqref{3.8}$_3$ over $(0,t)$ we get
\be
\bfI\cdot{\bfomega}_\infty(t)=\bfI\cdot{\bfomega}_\infty(0)-\int_0^t(\bfomega_\infty+\bfa)\times(\bfI\cdot\bfomega_\infty)\,,\ \ \mbox{for all $t\in(0,\infty)$}\,.
\eeq{4.2}
\begin{definition} A pair $(\bfv,\bfomega_\infty)$ is a {\em weak solution} to \Eqref{3.6}, \Eqref{3.8}--\Eqref{3.9} if it satisfies the following properties
\begin{itemize}
  \item[{\rm (i)}] $\bfv\in C_w([0,\infty);H(\mathscr C))\cap L^\infty(0,\infty;H(\mathscr C))\cap L^2(0,\infty;\cald_0^{1,2}(\mathscr C))$\,; 
  \item[{\rm (ii)}] $\bfomega_\infty\in C^{1}((0,\infty))\cap C([0,\infty))$\,; 
  \item [{\rm (iii)}] Strong energy inequality:\footnote{Notice that, in view of \Eqref{3.12}, every sufficiently regular solution satisfies \Eqref{4.3} with the equality sign (and for all $s\ge0$).}  
\be  
\mathcal E(t)+\mu\int_s^t\|\nabla\bfv(\tau)\|_2^2\le \mathcal E(s)
\eeq{4.3}
for $s=0$, for a.a. $s>0$ and all $t\ge s$,  where $\mathcal E$ is defined in \Eqref{3.13}\,;
\item[{\rm (iv)}] $(\bfv,\bfomega_\infty)$ obey \Eqref{4.1}--\Eqref{4.2}\,. 
\end{itemize}
\end{definition}

Combining the classical Galerkin method with the a priori estimate obtained after integration of \Eqref{3.12} in time, by more or less standard arguments one can show that the class of weak solutions is not empty.  In this regard, we have the following existence result for whose proof we refer the interested reader to \cite[Chapter 3]{Mazz}; see also \cite[Theorem 5.6]{ST}.
\Bt For any given $\bfv_0\in H(\mathscr C)$, $\bfomega_{\infty0}\in\real^3$ there is at least one corresponding weak solution such that
$$
\lim_{t\to 0^+}\left(\|\bfv(t)-\bfv_0\|_2+|\bfomega_\infty(t)-\bfomega_{\infty0}|\right)=0\,.
$$
\ET{1}
\Br
Stated in physical terms, \theoref{1} ensures the existence of a corresponding weak solution, provided only the initial data have a finite kinetic energy.
\ER{1}
\Br Let $(\bfv,\bfomega_\infty)$,  be any weak solution, and set
$$
S=\{s\in [0,\infty): \ \mbox{the strong energy inequality \Eqref{4.3} holds.}\}\,.
$$
Then 
for arbitrarily given $\eta,\varepsilon>0$ there is $t_0\in S$ such that
$$
\Int{t_0}{\infty}\|\nabla\bfv(\tau)\|_2^2<\eta\,,\ \ \|\nabla\bfv(t_0)\|_2<\varepsilon\,.
$$
In fact, from the strong energy inequality \Eqref{4.3} and the fact that $[0,\infty)$ and $S$ have the same one-dimensional Lebesgue measure, it follows that
\be
\int_{0}^\infty\|\nabla\bfv(\tau)\|_2^2=\int_S\|\nabla\bfv(\tau)\|_2^2<\infty\,.
\eeq{4.4}
From this we deduce that there is at least one increasing, unbounded sequence $\{t_k\}\subset S$ with the property that for any $\varepsilon >0$ there is $\bar{k}\in \nat$ such that 
\be
\|\nabla\bfv(t_k)\|_2<\varepsilon\,,\ \ \mbox{for all $k\ge \bar k$.}
\eeq{4.5} 
Moreover, again by \Eqref{4.4} we infer that for any $\eta>0$ there is $\bar{t}>0$ such that
\be
\int_{\bar t}^\infty\|\nabla\bfv(\tau)\|_2^2<\eta\,.
\eeq{4.6}
As a result, the claimed property  follows from \Eqref{4.5}--\Eqref{4.6} by choosing $t_0=t_{k^*}$ where $k^*\ge \bar{k}$,  $t_{k^*}\ge \bar t$.
\ER{2}

We conclude this section with a continuous data-dependence result. As in the classical Navier--Stokes case, also in the present situation continuous dependence  of weak solutions upon initial data (or even uniqueness) is an open question. However, again as in the classical case, one can show that the above property holds for any weak solution possessing  a (further) mild degree of regularity. More precisely, we have the following result for whose proof we refer to \cite[\S 3.4]{Mazz}.  
\Bt Let $(\bfv,\bfomega_\infty)$ and $(\bfv^*,\bfomega^*_\infty)$ be two weak solutions corresponding to initial data $(\bfv_0,\bfomega_{\infty0})$ and $(\bfv^*_0,\bfomega^*_{\infty0})$, respectively. Suppose there is $T>0$ such that
\be 
\bfv^*\in L^q(0,T;L^r(\mathscr C))\,,\ \ \frac{2}{q}+\frac 3r=1\,,\ \ \mbox{some $r>3$}\,. 
\eeq{4.7}
Then, necessarily 
$$\ba{rl}\smallskip
\|\bfv(t)-\bfv^*(t)\|_2+|\bfomega_\infty(t)-\bfomega^*_{\infty}(t)|\le c\,\big(\|\bfv_0&-\bfv^*_0\|_2+|\bfomega_{\infty0}-\bfomega^*_{\infty0}|\big)\,,\\ &\mbox{for all  $t\in [0,T]$}\,,
\ea$$
where $c$ depends on $\essup{t\in [0,T]}\|\bfv(t)\|_2$, $\essup{t\in [0,T]}\|\bfv^*(t)\|_2$, $\|\bfv^*\|_{L^q(0,T;L^r(\mathscr C))}$, and $\Max{t\in[0,T]}|\bfomega^*_{\infty}(t)|$. Thus, in particular, if $(\bfv,\bfomega_\infty)$ and $(\bfv^*,\bfomega^*)$ have the same initial data, it follows that $(\bfv,\bfomega_\infty)\equiv(\bfv^*,\bfomega^*_\infty)$ a.e. in $\mathscr C\times [0,T]$.
\ET{2}\renewcommand{\theequation}{5.\arabic{equation}}
\setcounter{equation}{0}
\section{Large-Time Properties of Weak Solutions. Preliminary Results}
We begin to show that, for any weak solution, the kinetic energy of the liquid must decay to zero. More specifically, we have the following.
\Bp Let $(\bfv,\bfomega_\infty)$ be any weak solution. Then
$$
\lim_{t\to\infty}\|\bfv(t)\|_2=0\,.
$$
Moreover, if in particular all eigenvalues of the inertia tensor $\bfI$ coincide, namely, $\bfI=\lambda\,\bf1$, some $\lambda>0$, then there are constants $c_i>0$, $i=1,2$, such that
\be
\|\bfv(t)\|_2\le c_1\|\bfv(0)\|_2\,{\rm e}^{- c_2\,\mu\,t}\,,\ \ \mbox{all $t>0$\,,}
\eeq{SFX}
\EP{1}
{\bf Proof.} From the strong energy inequality \Eqref{4.3} we derive
\be\ba{rl}\smallskip
E(t)+\mu\Int st\|\nabla\bfv(\tau)\|_2^2&\le E(s)-\half(\bfomega_\infty(t)\cdot\bfI\cdot\bfomega_\infty(t)-\bfomega_\infty(s)\cdot\bfI\cdot\bfomega_\infty(s))\,,
\\
&\equiv E(s)+F
\ea
\eeq{5.1}
where $E$ is defined in \Eqref{3.14}. 
By the Poincar\'e inequality, 
\be
\|\bfv\|_2\le c_0\,\|\nabla\bfv\|_2\,,
\eeq{5.2}
so that, employing \Eqref{3.14}, the inequality \Eqref{5.1} implies
\be
E(t)+c_0\mu\int_s^tE(\tau)\le E(s)+F\,.
\eeq{5.3}
Next, from \Eqref{4.2} and Definition 1(ii) we infer
\be
\frac{1}{2}\ode{}{t}(\bfomega_\infty\cdot\bfI\cdot\bfomega_\infty)=-\bfa\times(\bfI\cdot\bfomega_\infty)\cdot\bfomega_\infty\,.
\eeq{5.4}
As a consequence, if $\bfI=\lambda\,\bf1$, we have $F\equiv 0$, and the second statement of the proposition follows from \lemmref{1} and \Eqref{3.14}. In the general case, from \Eqref{5.4},  \Eqref{3.9} and \Eqref{5.2} we deduce
\be
|F|\le c\int_s^t|\bfomega_\infty(\tau)|^2\|\bfv(\tau)\|_2\le c_1\int_s^t|\bfomega_\infty(\tau)|^2\|\nabla\bfv(\tau)\|_2
\eeq{5.5}
However, from the strong energy inequality \Eqref{4.3} with $s=0$, we derive $|\bfomega_\infty(t)|\le M$, with $M$ depending only on the initial data, and physical and geometric properties of $\mathscr S$. Consequently, from \Eqref{5.3} and \Eqref{5.5} we conclude 
$$
E(t)+c_0\mu\int_s^tE(\tau)\le E(s)+c_2\int_s^t\|\nabla\bfv(\tau)\|_2\,,
$$
so that the first statement of the proposition follows from this inequality, \Eqref{4.3}, \Eqref{3.14}, and \lemmref{1}.
\QED

The next result shows that, for large times, any weak solution becomes regular and, moreover, the velocity field decays to zero in the $W^{1,2}$-norm.
\Bp
Let $(\bfv,\bfomega_\infty)$ be a weak solution. Then, there exists $t_0>0$, with $t_0\to\infty$ as $\nu:=\mu/\rho\to0$,  such that
\be\ba{ll}\smallskip
\bfv\in C([t_0,T];\cald_0^{1,2}(\mathscr C))\cap L^\infty(t_0,\infty;\cald_0^{1,2}(\mathscr C))\cap L^2(t_0,T;W^{2,2}(\mathscr C))\,,\\ \partial_t\bfv\in L^2(t_0,T;L^{2}(\mathscr C))\,,\ \ \bfomega_\infty\in W^{2,2}(t_0,T)\,,
\ea\eeq{class}
for all $T>t_0$. Moreover, there exists $p\in L^2(t_0,T;W^{1,2}(\mathscr C))$, all $T>t_0$, such that $(\bfv,p,\bfomega_\infty)$ satisfy \Eqref{3.8}$_{1,2}$, a.e. in $\mathscr C\times (0,\infty)$. Finally,
\be
\lim_{t\to\infty}\|\bfv(t)\|_{1,2}=0\,.
\eeq{5.6}
\EP{2}
{\bf Proof.} Denote by $(\bfv_0,\bfomega_{\infty0})$ the initial data of the weak solution. Moreover, let $t_0\in S$ be an instant of time such that $\|\nabla\bfv(t_0)\|_2<\infty$, where $S$ is defined in \remref{2}. By that remark, we also know that such a $t_0$ certainly exists. We then construct a  {\em local}, in principle, strong solution with initial data $(\bfv(t_0),\bfomega_\infty(t_0))$. As shown in \cite[Chapter 4]{Mazz} (see also \cite[Theorem 7.1]{ST}), this can be accomplished by using the Galerkin method in conjunction with suitable estimates on the solution. For the reader's sake, and also for their important role played successively, we shall (formally) re-derive such estimates and, in fact, in a  simpler fashion than in \cite{Mazz}.  
By dot-multiplying both sides of \Eqref{3.8}$_1$ by $\partial_t\bfv$, integrating over $\mathscr C$, and taking into account \Eqref{3.6} and \Eqref{3.8}$_2$ we deduce (with $\nu:=\mu/\rho$)
\be\ba{rl}\smallskip
\Frac{\nu}{2}\ode{}{t}\|\nabla\bfv\|_2^2+
\|\partial_t\bfv\|_2^2+\dot{\bfa}&\cdot\Int{\mathscr C}{}\bfx\times\partial_t\bfv=-(\bfv\cdot\nabla\bfv,\partial_t\bfv)\\
&-\dot{\bfomega}_\infty\cdot\Int{\mathscr C}{}\bfx\times\partial_t\bfv-2(\bfomega_\infty+\bfa)\cdot\Int{\mathscr C}{}\bfv\times\partial_t\bfv\,.
\ea
\eeq{5.7}
From \Eqref{3.9} and  \lemmref{3}, it follows that
\be\ba{rl}\smallskip
\|\partial_t\bfv\|_2^2+\dot{\bfa}\cdot\Int{\mathscr C}{}\bfx\times\partial_t\bfv&=\|\partial_t\bfv\|_2^2-\left(\rho\bfI^{-1}\cdot\Int{\mathscr C}{}\bfx\times\partial_t\bfv\right)\cdot\left(\Int{\mathscr C}{}\bfx\times\partial_t\bfv\right)\\
&\ge c_1\|\partial_t\bfv\|_2^2\,.
\ea\eeq{5.8}
As a consequence, using \Eqref{5.8} into \Eqref{5.7} along with several applications of the Cauchy-Schwartz inequality, we show that
\be
\nu\,\ode{}{t}\|\nabla\bfv\|_2^2+c_2\|\partial_t\bfv\|_2^2\le c_3\,\left[|\dot{\bfomega}_\infty|^2+(|\bfa|^2+|\bfomega_\infty|^2)\|\bfv\|_2^2+\|\bfv\cdot\nabla\bfv\|_2^2\right]\,,
\eeq{5.9}
where, unless otherwise stated, here and in the following, the  quantities $c,c_i$, $i\in\nat$, denote constants independent of $\nu$.
Likewise, applying the projection operator $\mathscr P$ to both sides of \Eqref{3.8}$_1$, and arguing as before, we get
\be
\nu\,\|\mathscr P\Delta\bfv\|_2^2\le c_4\left[\|\partial_t\bfv\|_2^2+|\dot{\bfomega}_\infty|^2+(|\bfa|^2+|\bfomega_\infty|^2)\|\bfv\|_2^2+\|\bfv\cdot\nabla\bfv\|_2^2\right]
\eeq{5.10}
Thus, multiplying either side of \Eqref{5.10} by $c_2/2c_4$ and then adding side by side the resulting inequality and \Eqref{5.9}, we derive
\be\ba{rl}\smallskip
\nu\,\ode{}{t}\|\nabla\bfv\|_2^2+c_5\,\big(\|\partial_t\bfv\|_2^2&+\nu\,\|\bfv\|_{2,2}^2\big)\\&\le c_6\left[|\dot{\bfomega}_\infty|^2+(|\bfa|^2+|\bfomega_\infty|^2)\|\bfv\|_2^2+\|\bfv\cdot\nabla\bfv\|_2^2\right]\,,\ea
\eeq{5.11}
where we have used the well-known inequality
(e.g. \cite[\S IV.6]{Gabook})
$$
\|\bfv\|_{2,2}\le c\,\|\mathscr P\Delta\bfv\|_2\,.
$$
Our next objective is to estimate the right-hand side of \Eqref{5.11}. By Sobolev embedding theorem and Cauchy-Schwartz inequality
\be
\|\bfv\cdot\nabla\bfv\|_2^2\le \|\bfv\|_\infty^2\|\nabla\bfv\|_2^2\le c_7\,\|\nabla\bfv\|_2^3\|\bfv\|_{2,2}\le \frac{c_8}{\nu}\|\nabla\bfv\|_2^6+\frac{c_5\,\nu}{2c_6}\|\bfv\|_{2,2}^2\,.
\eeq{5.12}
Furthermore, by the strong energy inequality \Eqref{4.3} with $s=t_0$, \Eqref{3.9}, and \Eqref{3.14} we readily show that
\be
|\bfa(t)|+|\bfomega(t)|+\|\bfv(t)\|_2\le c_9\,\big(|\bfomega_\infty(t_0)|+\|\bfv(t_0)\|_2\big)\,,
\eeq{5.13}
whereas \Eqref{3.8} and the preceding inequality imply
\be
|\dot{\bfomega}_\infty(t)|\le c_{10}\,\big(|\bfomega_\infty(t_0)|^2+\|\bfv(t_0)\|_2^2\big)\,. 
\eeq{5.14}
However, again by the strong energy inequality  \Eqref{4.3} with $s=0$ and $t=t_0$, \Eqref{3.9}, and \Eqref{3.14} we deduce
$$
|\bfomega_\infty(t_0)|+\|\bfv(t_0)\|_2\le c_{11}\,\big(|\bfomega_{\infty0}|+\|\bfv_0\|_2\big)\,,
$$
so that, from the latter and \Eqref{5.13}--\Eqref{5.14} we conclude
\be
|\dot{\bfomega}_\infty(t)|^2+(|\bfa(t)|^2+|\bfomega_\infty|^2)\|\bfv(t)\|_2^2\le k\,,
\eeq{5.15}
where $k$ depends only on the initial data of the weak solution, and  physical and geometric properties of $\mathscr S$.
Replacing \Eqref{5.12} and \Eqref{5.15} into \Eqref{5.11} produces 
\be
\nu\,\ode{}{t}\|\nabla\bfv\|_2^2+k_1\,\big(\|\partial_t\bfv\|_2^2+\nu\,\|\bfv\|_{2,2}^2\big)\le \frac{k_2}{\nu}\|\nabla\bfv\|_2^6+k_3\,,
\eeq{5.16}
where $k_i$, $i=1,2,3$, depend only on $\bfv_0,\bfomega_{\infty0}$, and physical and geometric properties of $\mathscr S$, but are independent of $\nu$. Inequality \Eqref{5.16} is the key tool that, as we mentioned earlier on, once combined with the classical Galerkin method produces the existence of a strong solution in a right neighborhood of $t_0$ \cite[Chapter 4]{Mazz}. More precisely, one shows that there exist $T^*>0$, 
and  a pair $(\tilde\bfv,\tilde{\bfomega}_\infty)$ such that
\be
\left.\ba{ll}\smallskip
\tilde\bfv\in C([t_0,t_0+\tau];\cald_0^{1,2}(\mathscr C))\cap L^\infty(t_0,t_0+\tau;\cald_0^{1,2}(\mathscr C))\\ \smallskip\tilde\bfv\in  L^2(t_0,t_0+\tau;W^{2,2}(\mathscr C))\,,\  \partial_t\tilde\bfv\in L^2(t_0,t_0+\tau;L^{2}(\mathscr C))\\ \tilde{\bfomega}_\infty\in C^2([t_0,t_0+\tau])
\ea\right\}\,\  \mbox{all $\tau\in(0,T^*)$.}\eeq{5.17}
Moreover, there is $p\in L^2(t_0,t_0+\tau;W^{1,2}(\mathscr C))$, all $\tau\in(0,T^*)$, such that $(\tilde\bfv,p,\tilde{\bfomega}_\infty)$ satisfy \Eqref{3.8}$_{1,2}$, a.e. in $\mathscr C\times (t_0,t_0+\tau)$. Let us continue to denote by $T^*$ the largest positive real-extended number for which \Eqref{5.17} holds. Then either $T^*=\infty$, or else, if $T^*<\infty$, necessarily 
\be
\lim_{t\to T_1^+}\|\nabla\tilde\bfv(t)\|_2=\infty\,,\ \ T_1:= t_0+T^*
\eeq{5.18}
The latter is shown by a classical argument. In fact, suppose there  is a sequence $\{t_m\}\subset [t_0,T_1)$ with $t_m\to T_1^+$ and such that
\be
\|\nabla\tilde\bfv(t_m)\|_2\le M\,,
\eeq{5.19}
where $M$ independent of $m$. Setting $z:=\|\nabla\tilde{\bfv}\|_2^2+k_3^{\frac{1}{3}}$, 
from \Eqref{5.16} one shows
$$
z'(t)\le c_{12}\,z^3(t)\,,
$$
with $c_{12}$ depending also on $\nu$, 
which in turn furnishes
$$
z^2(t)\le\frac{z^2(t_m)}{1-c_{12}z^2(t_m)(t-t_m)}\,.
$$
Using this inequality and \Eqref{5.19} it immediately follows
$$
z^2(t)\le c_{13}\,,\ \ t\in[t_m,t_m+M_1]\,,
$$
$M_1$ independent of $m$, which, by taking $m$ sufficiently large, proves $\|\nabla\bfv (t)\|\le c_{14}$ for all $t\in [t_m,T_2]$ with $T_2>T_1$. By employing the method previously described we can then extend the solution in the class \Eqref{5.17} to a time interval $[0,T_2]$, with $T_2>t_0+T^*$, which contradicts the assumption that $[t_0.t_0+T^*)$ is maximal with $T^*<\infty$. We next observe that, by the uniqueness result in  \theoref{2},  our original weak solution $(\bfv,\bfomega_\infty)$  must coincide with the strong solution just found, $(\tilde{\bfv},\tilde{\bfomega}_\infty)$ on the time interval $[t_0,t_0+T^*)$. We shall now show that, by choosing $t_0$ appropriately, \Eqref{5.18} (with $\tilde{\bfv}\equiv\bfv$) cannot occur, thus implying $T^*=\infty$, which completes the proof of the proposition.
To this end, we observe that from \Eqref{5.16} it follows, in particular, that
\be
\ode{}{t}\|\nabla\bfv\|_2^2\le -k_1\|\nabla\bfv\|_2^2 +\frac{k_2}{\nu^2}\|\nabla\bfv\|_2^6+\frac{k_3}{\nu}\,.
\eeq{5.20}
We now pick $\delta\in(0,\sqrt{k_1\,\nu^2/k_2})\equiv (0,\kappa\,\nu)$, and choose $t_0\in S$ in such a way that
\be
\int_{t_0}^\infty\|\nabla\bfv(\tau)\|_2^2<\frac{\nu\,\delta^2}{4k_3}\,,\ \ \|\nabla\bfv(t_0)\|_2^2<\frac{\delta}{\sqrt{2}}\,.
\eeq{tvb}
In view of \remref{2} such a $t_0$ definitely exists. Moreover, by the above choice of $\delta$ and \Eqref{tvb}$_1$, it is also clear that $t_0\to\infty$ as $\nu\to0$. We may now apply \lemmref{2} with $t_1=\infty$ to infer\footnote{Notice that since $\bfv$ is in the class \Eqref{5.17}, it is well known that the map $t\to\|\nabla\bfv(t)\|_2$ is absolutely continuous; see, e.g., \cite[Lemma 1]{Prodi}.}
$$
\|\nabla\bfv(t)\|_2^2<\delta\,,\ \ \mbox{for all $t\ge t_0$}
$$ 
which shows that \Eqref{5.18} (with $\tilde\bfv\equiv\bfv$) cannot occur. Therefore, we conclude $T^*=\infty$, and this ends the proof of the existence  part of the proposition.  Finally,  \Eqref{5.6} follows again from \lemmref{2}, since in our case we may take $t_1=\infty$. 
\QED 
\begin{remark} The previous theorem shows that the time, $t_c$, after which the relative velocity field of the liquid becomes ``small'', is larger the smaller the kinematic viscosity. In view of \theoref{3} shown in Section 7, this also implies that the coupled system $\mathscr S$ will basically rotate uniformly as a single rigid body after an instant of time possessing the same property. Numerical experiments reported in Section 9.1 confirm this trend and, in addition,  provide an explicit  correlation between $ t_c$ and $\nu$.  
\end{remark}

We conclude this section with a useful result, obtained by combining \theoref{2} with \propref{2}, and whose proof is an immediate consequence of the fact that if $\bfv$ is in the class defined in \Eqref{class}, then it is also in the class defined in \Eqref{4.7}.
\Bc Let $(\bfv,\bfomega_\infty)$ be a weak solution, and $t_1\ge t_0$, where $t_0$ is defined in \propref{2}. Then, if $(\bfV,\bfW)$ is another weak solution with $(\bfv(t_1),\bfomega_\infty(t_1))=(\bfV(t_1),\bfW(t_1))$, necessarily $(\bfv,\bfomega_\infty)\equiv(\bfV,\bfW)$ a.e. in $\mathscr C\times [t_1,\infty)$. 
\EC{1}\renewcommand{\theequation}{6.\arabic{equation}}
\setcounter{equation}{0}
\section{Existence of the $\Omega-$Limit Set}
In this section we continue our investigation
of the asymptotic properties of a weak solution. More specifically, we shall show  that the ``ultimate'' motion of the coupled system solid-liquid must be a rigid motion.

To this end, following the standard literature, we define the $\Omega-$limit set of a weak solution $(\bfv,\bfomega_\infty)$.
\begin{definition} Let $(\bfv,\bfomega_\infty)$ be a weak solution. The  corresponding $\Omega-$limit set, $\Omega(\bfv,\bfomega_\infty)$, is  constituted by pairs $(\bfV,\bfp)\in H(\Omega)\times\real^3$ for which there exists an increasing and unbounded sequence $\{t_n\}\subset (0,\infty)$ such that 
$$
\lim_{n\to\infty}\left(\|\bfv(t_n)-\bfV\|_2+|\bfomega_\infty(t_n)-\bfp|\right)=0\,. 
$$
\end{definition}
\Br Notice that since we do not know whether weak solutions are unique, the $\Omega-$limit set is not just characterized by the initial data but, rather, by the particular weak solution that we are considering.
\ER{3}

In order to provide a first characterization of the $\Omega-$limit set, we  need the following simple but useful preparatory result.
\Bl
Let $\bfu \in L^\infty(0,\infty)$ be uniformly Lipschitz-continuous.
Then, given an unbounded increasing sequence $\{t_k\}\subset [0,\infty)$, there is  a subsequence $\{t_{k'}\}$ and a uniformly Lip\-sch\-itz-continuous function $\bar \bfu$ such that
\be 
\lim_{k'\to\infty}\bfu(t+t_{k'})=\bar \bfu(t)\,, \ \ \mbox{uniformly in $t\ge0$.}
\eeq{6.1}
\EL{4}{\bf Proof.} Let $\bfu_k(t):=\bfu(t+t_k)$, $k\in\nat$. Then, by assumption, $\{\bfu_k\}$ is equi-bounded. Moreover, by the uniform Lipschitz continuity property of $\bfu$, we get
\be 
|\bfu_k(t)-\bfu_k(s)|=|\bfu(t+t_k)-\bfu(s+t_k)|\le c\,|t-s|\,.
\eeq{6.2}
Thus the existence of a subsequence $\{t_{k'}\}$, and a continuous function $\bar\bfu$ satisfying \Eqref{6.1}  follows from Ascoli--Arzel\`a theorem. Finally, the uniform Lipschitz-continuity of $\bar\bfu$ is a consequence of \Eqref{6.1} and \Eqref{6.2}. 
\QED

We our now in a position to prove  the following preliminary properties of the $\Omega-$limit set of a generic weak solution.
\Bp Let $(\bfv,\bfomega_\infty)$ be a weak solution. Then
\be
\Omega(\bfv,\bfomega_\infty)=\{\0\}\times \cala\,,
\eeq{6.3}
where $\cala$ is a compact, connected subset of $\real^3$ that is left invariant by the semigroup relative to the system of differential equations
\be
\bfI\cdot\dot{\bfw}+\bfw\times(\bfI\cdot\bfw)=\0\,,\ \ \bfw\in C^1((0,\infty))\cap C([0,\infty))\,.
\eeq{6.4}
\EP{3}
{\bf Proof.} In view of the strong energy inequality \Eqref{4.3} --ensuring the uniform boundedness of $\bfomega_\infty(t)$--  and  \propref{1}, it is immediately proved that $\Omega(\bfv,\bfomega_\infty)$ is not empty, and  of the type given in \Eqref{6.3} where  $\cala$ is the compact set of $\real^3$ defined by the property that $\bfp\in\cala$ if and only if $\bfomega_\infty(t_k)\to \bfp$ for some increasing and unbounded $\{t_k\}\subset(0.\infty)$. The fact that $\cala$ is connected follows from standard arguments that, however, for the reader's sake we shall briefly reproduce here. If $\cala$  were not connected, since it is closed, we should have $\cala=S_1\cup S_2$, with  $S_1$ and $S_2$ closed and disjoint sets.  Therefore, setting $\varepsilon:=\dist(S_1,S_2)$, there would be two  sequences $\{s_k\}$ and $\{t_k\}$ such that for all sufficiently large $k$,  $s_k<t_k<s_{k+1}$ and
$$
\dist(\bfomega_\infty(s_k),S_1)<\frac\varepsilon 4\,,\ \ \dist(\bfomega_\infty(t_k),S_2)<\frac\varepsilon 4\,.
$$
Being $t\to\bfomega_\infty(t)$  continuous, we may find $\tau_k\in (s_k,t_k)$ such that
$$
\dist(\bfomega_\infty(\tau_k),S_1\cup S_2)\ge\frac\varepsilon 2\,.
$$
However, since $\bfomega(\tau_k)$ is bounded, by the latter displayed relation we can find $\bfpsi\not\in\cala$ and   unbounded  $\{\tau_{k'}\}$ such that, $\bfomega(\tau_{k'})\to \bfpsi$, which is at odds with the definition of $\cala$.  
It remains to show the stated invariance property. Let $\bfp\in\cala$, and denote by $\bfw=\bfw(t)$
the unique solution to \Eqref{6.4} with $\bfw(0)=\bfp$. We have to show that, at each time $t\ge 0$, $\bfw(t)\in \cala$. To this end it is enough to prove that there exists an increasing, unbounded sequence $\{t_n\}$ such that $\bfw(t)=\lim_{n\to\infty} \bfomega_\infty(t+t_n)$, for all the above values of $t$. Let $\{t_n\}$ be the sequence such that
\be
\bfp=\lim_{n\to\infty }\bfomega_\infty(t_n)\,.
\eeq{6.5}
From \Eqref{4.2} we then deduce
\be
\bfI\cdot \bfomega_\infty(t+t_n)=\bfI\cdot \bfomega_\infty(t_n)-\int_0^{t}\big\{[\bfa(\tau+t_n)+\bfomega_\infty(\tau+t_n)]\times\bfI\cdot \bfomega_\infty(\tau+t_n)\big\}\,d\tau\,.
\eeq{6.6}
Taking into account that, again by \Eqref{4.2} and \Eqref{4.3}, $\bfomega_\infty$ satisfies the assumptions of  \lemmref{4}, and that by \propref{1} $\bfa(t)\to 0$ as $t\to\infty$,  we may pass to the limit $n\to\infty$ in \Eqref{6.6}, possibly along a subsequence still denoted by $\{t_n\}$,   and use \lemmref{4} in conjunction with \Eqref{6.5} to infer
\be
\bfI\cdot \bar \bfu(t)=\bfI\cdot \bfp+\int_0^{t}\bar \bfu(\tau)\times \bfI\cdot \bar \bfu(\tau)\,d\tau\,,
\eeq{6.7}
for some $\bar \bfu\in C^{1}((0,\infty))\cap C([0,\infty))$, with $\bar{\bfu}(t)=\lim_{n\to\infty} \bfomega_\infty(t+t_n)$, all $t\in (0,\infty)$. 
Therefore, $\bfw$ and $\bar\bfu$ satisfy the same equation with the same initial condition, so that  by uniqueness we conclude $\bfw\equiv \bar \bfu$, and the proof is completed.
\QED
\renewcommand{\theequation}{7.\arabic{equation}}
\setcounter{equation}{0}
\section{Invariance of the $\Omega-$Limit Set and Main Result}
\propref{3} ensures that,  asymptotically,  the coupled system solid-liquid $\mathscr S$  will move  about its center of mass like a whole rigid body. The main objective of this section is to refine this result, and  prove that, in fact, such a rigid motion must reduce to a permanent rotation around one of the central axes of inertia of $\mathscr S$.

A crucial point to derive the above conclusion is to secure that the $\Omega-$limit set is invariant {\em in the class of weak solutions}. In view of \propref{3}, the latter amounts to say that 
every weak solution $(\bfV,\bfW)$ (say) corresponding to initial data $(\0,\bfp)\in \Omega(\bfv,\bfomega_\infty)$ satisfies $(\bfV(t)\equiv\0,\bfW(t))\in \Omega(\bfv,\bfomega_\infty)$, for all $t>0$.

As is well known, the invariance property, typically, requires the uniqueness of the solution, a feature  that, in the present case, is not available.~\footnote{In the specific case of the Navier--Stokes equations, the uniqueness request can be  
relaxed to an a priori weaker condition like, for example, continuity in the ``energy'' norm \cite{Ball}, which, however, it is still an unproved property for weak solutions.} 
Nevertheless, using the fact that the velocity field of the liquid decays asymptotically to zero, we shall prove that  to infer the desired invariance, a sort of ``asymptotic uniqueness'' is indeed  enough. In turn, the latter will follow from  the fact that, by \propref{2}, every weak solution becomes regular for sufficiently large times. 

In order to accomplish all the above, we begin to define the following set
\be
\calt:=\{\bfw\in\real^3:\ \bfw\times(\bfI\cdot\bfw)=\0\}\,.
\eeq{7.1}
We recall that a motion of the coupled system solid-liquid $\mathscr S$ is called a {\em permanent rotation} if $\mathscr S$ rotates as a whole rigid body with a constant angular velocity. Then, 
from the physical viewpoint, $\calt$ characterizes the class of all admissible angular velocities, $\bfw$, in a permanent rotation of $\mathscr S$. Notice that from \Eqref{7.1} it follows that $\bfw$ must be directed along an eigendirection of the tensor $\bfI$. Thus,  if the three eigenvalues $\lambda_i$, $i=1,2,3$, of $\bfI$ are all distinct, then 
$$
\calt:=\{\alpha\,\bfe_1\}\cup\{\beta\,\bfe_2\}\cup\{\gamma\,\bfe_3\}\,,\ \ \alpha,\beta,\gamma\in\real\,,
$$     
where, we recollect, $\bfe_i$ is the (normalized) eigenvector corresponding to $\lambda_i$, $i=1,2,3$. Furthermore, if $\lambda_1=\lambda_2\neq \lambda_3$ (say), then 
$$ 
\calt=\{\bfz\in\real^2: \bfz=\alpha\,\bfe_1+\beta\,\bfe_2\,,\ \alpha,\beta\in\real\}\cup \{\gamma\,\bfe_3\}\,,\ \ \gamma\in\real.
$$
Finally, if $\lambda_1=\lambda_2=\lambda_3$, then $\calt=\real^3$.\smallskip\par

The following lemma ensures that the invariance of the $\Omega-$limit set in the class  of weak solutions reduces the set $\cala$ in \propref{3} to an element of the set $\calt$.
\Bl Let $(\bfv,\bfomega_\infty)$ be a given weak solution, and suppose that $\Omega(\bfv,\bfomega_\infty)$ is invariant in the class of weak solutions. Then
$$
\Omega(\bfv,\bfomega_\infty)=\{\0\}\times \{\bfw\}\,,\ \ \mbox{for some $\bfw\in\calt$.}
$$
\EL{5}
{\bf Proof.} Under the given assumptions, and in view of \propref{3}, the velocity field $\bfv$ of the generic weak solution with initial data in the $\Omega-$limit set must satisfy $\bfv(t)\equiv \0$. As a result, from \Eqref{4.1}--\Eqref{4.2} we derive
\be\ba{cc}\medskip
\Int{\mathscr C}{}\dot{\bfomega}_\infty\times \bfx\cdot\bfpsi=0\,,\ \ \mbox{for all $\bfpsi\in \cald_0^{1,2}(\mathscr C)$\,,}\\
\bfI\cdot\dot{\bfomega}_\infty+\bfomega_\infty\times(\bfI\cdot\bfomega_\infty)=\0\,.
\ea
\eeq{7.2}
From the first relation in \Eqref{7.2} it easily follows that
\be\dot{\bfomega}_\infty\times \bfx=\nabla \phi\,,
\eeq{7.3} 
with $\phi$ a suitable smooth scalar field. Thus, operating with $\nabla\times$ on both sides of \Eqref{7.3}, we deduce  that there is $\bfw\in\real^3$ such that $\bfomega_\infty(t)=\bfw$,  for all $t\ge 0$. Replacing this information back into \Eqref{7.2}$_2$ we conclude $\bfw\in\calt$, which completes the proof of the lemma. 
\QED

\smallskip\par
The next result proves that the assumption of the previous lemma is indeed met.
\Bp $\Omega(\bfv,\bfomega_\infty)$ is invariant in the class of weak solutions.
\EP{4}
{\bf Proof.} We begin to introduce some further notation. For $(\bar\bfv,\bar\bfomega)\in H(\mathscr C)\times \real^3$, we denote by
$$
(\bfu(t;\bar\bfv,\bar\bfomega), \bfq(t;\bar\bfv,\bar\bfomega))
$$
a weak solution\footnote{In principle there is more than one, due to lack of a uniqueness result.} at time $t\ge0$, with initial data $\bfu(0)=\bar\bfv$, $\bfq(0)=\bar\bfomega$. Moreover, if $\bfp\in\real^3$, we indicate by $\bfw(t;\bfp)$ the solution to \Eqref{6.4} at time $t\ge 0$,  with $\bfw(0)=\bfp$. Also, if $(\bfv_1,\bfomega_1), (\bfv_2,\bfomega_2)\in H(\mathscr C)\times \real^3$, we set
$$
\|(\bfv_1,\bfomega_1)-(\bfv_2,\bfomega_2)\|:=\|\bfv_1-\bfv_2\|_2+|\bfomega_1-\bfomega_2|\,.
$$
Finally, we indicate by $(\bfv_0,\bfomega_{\infty0})$ the initial data of the given weak solution $(\bfv,\bfomega_\infty)$.
Now, let $\bfp\in \mathcal A$. By \propref{3} we know that $\bfw(t;\bfp)\in\mathcal A$, for all $t\ge0$. The desired property  will then be accomplished if we show 
\be
(\bfV(t;\0,\bfp),\bfW(t;\0,\bfp))=(\0,\bfw(t;\bfp))\,,\ \ \mbox{for a.a. $t\ge 0$\,.}
\eeq{7.4}
Since $\bfp\in\cala$, there is $\{t_n\}\subset (0,\infty)$ increasing and unbounded such that
\be
\lim_{n\to\infty}|\bfomega_\infty(t_n)-\bfp|=0\,.
\eeq{7.5}
Again by \propref{3} we also deduce that 
$$
\lim_{n\to\infty}|\bfomega_\infty(t+t_n)-\bfw(t;\bfp)|=0\,,\ \ \mbox{for all $t\ge 0$.}
$$
Furthermore, from \propref{2} we infer
\be
\lim_{n\to\infty}\|\bfv(t+t_n)\|_2=0\,,\ \ \mbox{all $t\ge 0$.}
\eeq{7.6}
As a result, we conclude
\be
\lim_{n\to\infty}\|(\bfv(t+t_n;\bfv_0,\bfomega_{\infty0}),\bfomega_\infty(t+t_n;\bfv_0,\bfomega_{\infty0})-(\0,\bfw(t;\bfp))\|=0\,,\ \mbox{all $t\ge 0$}\,.
\eeq{7.7}
We next observe that, by \cororef{1}, for all sufficiently large $n$
$$\ba{rl}\medskip
(\bfv(t+t_n;\bfv_0,&\bfomega_{\infty0}),\bfomega_\infty(t+t_n;\bfv_0,\bfomega_{\infty0}))\\&=
(\bfv(t;\bfv(t_n),\bfomega_\infty(t_n)),\bfomega_\infty(t;\bfv(t_n),\bfomega_\infty(t_n))\,,\ \ \mbox{all $t\ge 0$}.
\ea
$$
Therefore, by \theoref{2},
\be\ba{rl}\medskip
\|(\bfV(t;\0,\bfp),\bfW(t;\0,\bfp))&-(\bfv(t+t_n;\bfv_0,\bfomega_{\infty0}),\bfomega_\infty(t+t_n;\bfv_0,\bfomega_{\infty0}))\|\\
&\le c\,\big(\|\bfv(t_n)\|_2+|\bfomega_\infty(t_n)-\bfp|\big)\,,\ \ \mbox{all $t\ge 0$.}
\ea
\eeq{7.8}
Since
$$\ba{rl}\medskip
\|(\bfV&(t;\0,\bfp),\bfW(t;\0,\bfp))-(\0,\bfw(t;\bfp))\|\\ \medskip\le& \|(\bfV(t;\0,\bfp),\bfW(t;\0,\bfp))-(\bfv(t+t_n;\bfv_0,\bfomega_{\infty0}),\bfomega_\infty(t+t_n;\bfv_0,\bfomega_{\infty0}))\|\\
&+\|(\bfv(t+t_n;\bfv_0,\bfomega_{\infty0}),\bfomega_\infty(t+t_n;\bfv_0,\bfomega_{\infty0}))-
(\0,\bfw(t;\bfp))\|\,,
\ea
$$
by letting $n\to\infty$ in the latter relation and employing  \Eqref{7.5}--\Eqref{7.8}, we show \Eqref{7.4}, which concludes the proof of the proposition.
\QED

We are now in a position to provide a complete description of the asymptotic behavior in time of the coupled system solid-liquid.
\Bt Let $\mathscr S$ be the coupled system constituted by a rigid body with an interior cavity $\mathscr C$ of class $C^2$  completely filled with a Navier--Stokes liquid. Suppose that no external forces act on $\mathscr S$. \par Let  $(\bfv,\bfomega_\infty)$ be any weak solution, in the sense of\, {\rm Definition 1}, to the initial-boundary value problem \Eqref{3.6}, \Eqref{3.8}--\Eqref{3.9} governing the motion of $\mathscr S$. Then, 
\be
\lim_{t\to\infty}\|\bfv(t)\|_{1,2}=0\,,
\eeq{7.9}
and there exists $\bar\bfomega\in\real^3$ such that
\be
\lim_{t\to\infty}\bfomega_\infty(t)=\bar{\bfomega}\,.
\eeq{7.10}
If, in particular, all eigenvalues of $\bfI$ coincide, then
\be\ba{rl}\medskip
\|\bfv(t)\|_{2}+|\bfomega_\infty(t)-\bar\bfomega|&\le c_1\,{\rm e}^{-c_2t}\,,\  \mbox{ all $t>0$\,,}\\
\|\nabla\bfv(t)\|_{2}&\le c_1\,{\rm e}^{-c_2t}\,,\  \mbox{ all sufficiently large $t>0$\,,}
\ea
\eeq{7.11}
for some $c_1,c_2>0$.

Furthermore, the vector $\bar\bfomega$, is parallel to one of the eigenvectors, $\bfe$, of the central inertia tensor $\bfI$. Finally, 
\be
\bar{\bfomega}=\frac{1}{\lambda}\,\bfK_G
\eeq{7.12}
where $\lambda$ is the eigenvalue of $\bfI$ associated with $\bfe$, representing the moment of inertia of $\mathscr S$ with respect to $\bfe$, and $\bfK_G$ is the (constant) angular momentum of $\mathscr S$ with respect to $G$.

Therefore, the asymptotic motion of $\mathscr S$ is 
a constant rigid rotation
around a central axis of inertia of $\mathscr S$ that aligns with the direction of the constant total angular momentum. 
\ET{3}
{\bf Proof.} Property \Eqref{7.9} is proved in \propref{2}. Moreover, by the strong energy inequality \Eqref{4.3} with $s=0$,  and \Eqref{3.14} we deduce that $|\bfomega_{\infty}(t)|$ is uniformly bounded in time. As a result, from every increasing, unbounded sequence $\{t_n\}\subset(0,\infty)$ we can select a subsequence (still denoted by $\{t_n\}$) such that
$$
\lim_{n\to\infty}\bfomega_{\infty}(t_n)\ \ \ \mbox{exists}.  
$$
However, from \propref{4} and \lemmref{5} we know that this limit is independent of the particular sequence $\{t_n\}$ and coincides with some vector $\bar \bfomega\in\calt$, with $\calt$ defined in \Eqref{7.1}. This proves \Eqref{7.10} and the first property of the vector $\bar\bfomega$.  As for \Eqref{7.12}, we observe that, by the conservation of angular momentum,
$$
\bfI\cdot\bfomega_\infty(t)=\sum_{i=1}^3\lambda_i\omega_{\infty i}(t)\bfe_i(t)=\bfK_G\,,\ \mbox{all $t\ge 0$}\,,
$$
where, we recall, $\lambda_i$ is an eigenvalue of $\bfI$ and $\bfe_i$ corresponding eigenvector, $i=1,2,3$.  
Consequently, by passing to the limit $t\to\infty$ in the latter relation and taking into account \Eqref{7.10} and that $\bar\bfomega\in\calt$, we show the validity of \Eqref{7.12}. It remains to prove the exponential decay, under the  assumption that $\bfI=\lambda\bf1$, $\lambda>0$. To this end, we notice that, by Schwartz inequality and \propref{1},
\be
|\bfa(t)|+\|\bfv(t)\|_2\le c_1\,{\rm e}^{-c_2t}\,.
\eeq{7.13} 
Further, from \Eqref{4.2}  we get
\be
\bfomega_\infty(t)-\bfomega_\infty(s)=\int_s^t\bfa(\tau)\times\bfomega_\infty(\tau)\,,\ \ \mbox{for all $t\ge s\ge0$.}
\eeq{7.14} 
Therefore, recalling that $|\bfomega_\infty(t)|$ is uniformly bounded, from \Eqref{7.10}, \Eqref{7.13}, and \Eqref{7.14} we show
\be
|\bfomega_\infty(t)-\bar\bfomega|\le c_3\int_t^\infty{\rm e}^{-c_2\tau}\le c_4 {\rm e}^{-c_2t}\,,\ \ \mbox{all $t>0$}\,.
\eeq{7.15} 
By \Eqref{7.14} and \Eqref{7.13} we also obtain
\be
|\dot{\bfomega}_\infty(t)|\le c_5\,{\rm e}^{-c_2t}\,,\ \ \mbox{all $t>0$}\,.
\eeq{7.16}
Employing \Eqref{7.13} and \Eqref{7.15}--\Eqref{7.16} into \Eqref{5.11} we then deduce, for all sufficiently large $t$,
\be
\ode{}{t}\|\nabla\bfv\|_2^2+c_6\|\bfv\|_{2,2}^2\le c_7\,\big({\rm e}^{-2c_2t}+\|\bfv\cdot\nabla\bfv\|_2^2\big)\,.
\eeq{7.17}
Next, by H\"older and Sobolev inequalities,
$$
\|\bfv\cdot\nabla\bfv\|_2\le c_8\,\|\bfv\|_2^{\frac14}\|\nabla\bfv\|_2^{\frac34}\|\bfv\|_{1,2}^{\frac14}\|\bfv\|_{2,2}^{\frac34}\,,
$$
so that using \Eqref{5.2},  and recalling that by \propref{4}, $\|\nabla\bfv(t)\|_2$ is uniformly bounded for sufficiently large $t$ we infer for all such times
$$
\|\bfv\cdot\nabla\bfv\|^2_2\le c_{9}\|\bfv\|_2^{\frac12}\|\bfv\|_{2,2}^{\frac32}\le c_{10}\,\|\bfv\|_2^2+\frac{c_6}{2}\|\bfv\|_{2,2}^2\,,
$$
where, in the last step, we made use of the Young inequality. From the latter relation, \Eqref{7.13}, and \Eqref{7.17} we thus derive, in particular,
\be
\ode{}{t}\|\nabla\bfv\|_2^2\le c_{11}\,{\rm e}^{-2c_2t}\,.
\eeq{7.18}
Integrating \Eqref{7.18} over $(t,\infty)$ and taking into account \Eqref{7.9} we show
$$
\|\nabla\bfv(t)\|_2^2\le c_{12}\,{\rm e}^{-2c_2t}\,,\ \mbox{all sufficiently large $t>0$}\,,
$$
which, once combined with \Eqref{7.13} and \Eqref{7.15} completes the proof of the theorem.
\QED
\Br
The previous theorem gives a full rigorous proof of Zhukhovskii's conjecture in a very general class of solutions, and for sufficiently smooth cavities.  
\ER{Z}
\begin{remark} From \Eqref{7.12} it follows that $\bar\bfomega=\0$ if and only if $\bfK_G=\0$. Notice that the latter condition is not physically relevant. Actually,   it is satisfied either by identically vanishing initial data, in which case the rest is the only corresponding weak solution, or else, more generally, for initial data  able to produce, at time $t=0$ an  angular momentum of the liquid (relative to the rigid body) that is {\em exactly} the opposite of that of the rigid body, a circumstance that is very unlikely to happen. We also observe that, from Definition 1 it immediately follows that if $\bfK_G=\0$ {\em every} weak solution must have $\bfomega_\infty(t)=\0$ for all $t\ge 0$. With the help of \Eqref{5.4}, \Eqref{3.14}, and \lemmref{1} this implies, in turn,  $\|\bfv(t)\|_2\le c_1 \|\bfv(0)\|_2\,{\rm e}^{-c_2 t}$, for some $c_1,c_2>0$ and all $t\ge 0$, thus re-obtaining, in a simpler way, the result of \cite[Theorem 5.6]{ST}.  
\label{rem:ST}
\end{remark}
\renewcommand{\theequation}{8.\arabic{equation}}
\setcounter{equation}{0}
\section{On the Attainability and Stability of Permanent Rotations}
As we have shown in the previous section, the system $\mathscr S$ will eventually perform a permanent rotation, as a single rigid body, around one of the central axes of inertia. However, our result does not specify around which axis this rotation will be attained. This issue assumes even more significance if we  keep in mind that weak solutions may lack of uniqueness and therefore, in principle, we may have two different solutions with the same initial data generating, asymptotically, two permanent rotations around different axes.  One of the objectives of this section is therefore to analyze this problem in some details. In particular, we shall prove that, if the initial data satisfy certain sufficiently general conditions, the   permanent rotation will always occur along that central axis with the largest moment of inertia; see \theoref{4}.  

The other related objective concerns the stability of such permanent rotations. In this respect, we recall that this problem has a long history, beginning with  the work of Sobolev \cite{Sob} (in the inviscid case),  Rumyantsev \cite[p. 203]{Rum} (see also \cite[\S 3-3]{MoRu}), Chernousko \cite{Ch}, and Smirnova \cite{Smirn1,Smirn2} until the more recent work of Kostyuchenko {\em et al.} \cite{KSY}, and Kopachevsky \& Krein \cite{KK}. However, besides \cite{Rum,MoRu}, 
results proved by these authors are obtained under different simplifications. More precisely, in \cite{KSY,KK}  stability/instability is studied by linearizing the relevant equations.  Similar results are obtained  in \cite{Ch,Smirn1} under the assumption of  small Reynolds number, which allows one to approximate  the problem with an appropriate system of nonlinear ordinary differential equations.  On the other hand, in  \cite{Rum} a nonlinear analysis is performed that provides sufficient conditions for stability of permanent rotations.  Our other main objective in this section is to show {\em necessary and sufficient} conditions for stability for the full nonlinear problem,  without any approximation. These conditions contain those of  \cite{Rum} as a particular case, and extend those of \cite{KSY,Ch,Smirn1} to the nonlinear level.  

In order to show all the above, let us begin to rename, as more customary,  the eigenvalues of $\bfI$ by $A$, $B$ and $C$, with corresponding eigenvectors $\bfe_1$, $\bfe_2$ and $\bfe_3$. We recall that, from the physical viewpoint, $A$, $B$ and $C$ are the moment of inertia of $\mathscr S$ around the axes passing through $G$ and parallel to $\bfe_1$, $\bfe_2$, and $\bfe_3$, respectively ({\em central moment of inertia}). Moreover,  set
$$
\bfomega_\infty=p\,\bfe_1+q\,\bfe_2+r\,\bfe_3\,.
$$

Our approach to attainability and stability is quite straightforward and  relies upon the following three ingredients: (i) \theoref{3}, (ii) balance of energy, and (iii) conservation of angular momentum. To this end, we begin to observe that the strong energy inequality \Eqref{4.3}  can be written as follows
\be\ba{rl}\smallskip
2E(t)+
A\,p^2(t)+&\!
B\,q^2(t)+C\,r^2(t)+2\mu\Int0t\|\nabla\bfv(\tau)\|_2^2\\ &
\le 2E(0)+A\,p^2(0)+B\,q^2(0)+C\,r^2(0)\,, \ \ \mbox{all $t\ge 0$}\,,\ea
\eeq{8.1}
where the ``energy'' $E$ is defined in \Eqref{3.14}. 
Furthermore, by dot-multiplying both sides of \Eqref{3.8}$_3$ by $\bfI\cdot\bfomega_\infty$ we obtain the following equation representing conservation of (the magnitude of) angular momentum
\be
A^2p^2(t)+B^2q^2(t)+C^2r^2(t)=A^2p^2(0)+B^2q^2(0)+C^2r^2(0)\,.
\eeq{8.2}

The next result concerns the attainability of permanent rotations. Without loss of generality, we shall assume throughout $A\le B\le C$. 
\Bt The following statements hold.\footnote{We assume  $\bfomega_\infty(0)\neq\0$, otherwise the motion of the coupled system is physically irrelevant; see Remark \ref{rem:ST}.}
\begin{itemize}
\item[{\rm (a)}] Suppose $A=B<C$. Then, if
\be
0<E(0)\le\frac{(C-A)\,C}{2 A}r^2(0)\,,
\eeq{8.3}
necessarily
\be\ba{ll}\medskip
\Lim{t\to\infty}p(t)=\Lim{t\to\infty}q(t)=0\\
\Lim{t\to\infty}r(t)=\bar{r}\neq 0\,.
\ea
\eeq{8.4}
\item[{\rm (b)}] Suppose $A<B=C$. Then, if 
\be
0<E(0)\le \frac{B(B-A)}{2 A}(q^2(0)+r^2(0))\,,
\eeq{8.5}
necessarily
\be\ba{ll}\medskip
\Lim{t\to\infty}p(t)=0\,,\ \ \Lim{t\to\infty}q(t)=\bar q\\
\Lim{t\to\infty}r(t)=\bar{r}\,,
\ea
\eeq{8.6}
where at least one of $\bar q,\bar r$ is not zero.
\item[{\rm (c)}] Suppose $A<B<C$. Then, if  
\be\ba{ll}\medskip
0<E(0)+\Frac{A}{2B}(B-A)p^2(0)\le \Frac C{2B}(C-B)r^2(0)\,,\\
0<E(0)\le\Frac{B}{2A}(B-A)q^2(0)+\Frac C{2A}(C-A)r^2(0)\,,
\ea
\eeq{PN2}
necessarily \Eqref{8.4} follows.
\end{itemize}
\ET{4}
{\bf Proof.}  We commence by noticing that from \theoref{3} we know that
\be
\lim_{t\to\infty}\bfomega_\infty(t)=\bar{p}\bfe_1+\bar{q}\bfe_2+\bar{r}\bfe_3\,,\ \ \lim_{t\to\infty}E(t)=0\,,
\eeq{PN4}
for some $\bar p,\bar q,\bar r\in\real$.
Thus, passing to the limit $t\to\infty$ on both sides of \Eqref{8.1} and \Eqref{8.2} we deduce
\be\ba{rl}\medskip
A\bar p^2+B\bar q^2+C\bar r^2+2\mu\Int0\infty\|\nabla\bfv(t)\|_2^2&\le 2E(0)+Ap^2(0)+Bq^2(0)+Cr^2(0)
\\
A^2\bar p^2+B^2\bar q^2+C^2\bar r^2
&=A^2p^2(0)+B^2q^2(0)+C^2r^2(0)
\,.
\ea
\eeq{PN5}In order to show (a), we observe that again
by \theoref{3} either   $\bar p=\bar q=0$ or $\bar r=0$. Let us show  that the latter cannot occur. In fact,
multiplying both sides of \Eqref{PN5}$_1$ by $A$ (=$B$),  subtracting  \Eqref{PN5}$_2$, side by side, to the resulting inequality  and taking $\bar r=0$, we deduce
$$
2A\mu\int_0^\infty\|\nabla\bfv(\tau)\|_2^2\le 2AE(0)+C(A-C)r^2(0)\,,
$$
which  cannot hold under the assumption \Eqref{8.3}. To show statement (b), suppose $\bar q=\bar r=0$. We then multiply  both sides of  \Eqref{PN5}$_1$ by $A$ and subtract to the resulting inequality equation \Eqref{PN5}$_2$, side by side, to get
\be
2A\mu\Int0\infty\|\nabla\bfv(t)\|_2^2\le 2AE(0)+B(A-B)q^2(0)+C(A-C)r^2(0)\,,
\eeq{SF1}
which is at odds with \Eqref{8.5}. It remains to demonstrate (c), namely, $\bar p=\bar q=0$. Suppose $\bar p\neq 0$. Then by \theoref{3}, $\bar q=\bar r=0$, and, therefore, \Eqref{SF1} holds, which is contradicted by \Eqref{PN2}$_2$.  
Suppose, instead, $\bar q\neq 0$. 
Then, again by \theoref{3}, $\bar p=\bar r=0$. Thus, multiplying  both sides of  \Eqref{PN5}$_1$ by $B$ and subtracting to the resulting inequality equation \Eqref{PN5}$_2$, side by side, we  infer 
\be
2B\mu\Int0\infty\|\nabla\bfv(t)\|_2^2\le 2BE(0)+A(B-A)p^2(0)+C(B-C)r^2(0)\,.
\eeq{SF2}
However,  \Eqref{SF2} is in contrast with \Eqref{PN2}$_1$, and  
the proof of the theorem is completed.
\QED

\Br

As we mentioned earlier on, the above theorem assumes great  relevance  when the coupled system $\mathscr S=\mathscr B\cup \mathscr L$ has gyroscopic structure around the $\bfe_3$ axis (say), that is, $A=B\neq C$. In such a case, our result ensures, in particular, that  if the liquid is initially at rest with respect to the rigid body (that is, the relative velocity field of the liquid is zero at $t=0$), eventually, the permanent rotation of $\mathscr S$ will occur along the axis of the gyroscope, $\{G,\bfe_3\}:={\sf a}$, {\em if and only if} the moment of  inertia with respect to that axis is larger than those around the other two. As an illustration of this fact, consider the case where the body $\mathscr B$ is a  hollow cylinder  (like a  metal can), completely filled with a viscous liquid $\mathscr L$. In this situation, {\sf a} coincides with the axis of the cylinder. We assume that the central moment of inertia of $\mathscr B$ are negligible compared to those of $\mathscr L$. Let $h$ and $R$  be height and radius of $\mathscr B$, respectively. Then combining \theoref{3} and \theoref{4}, we may state that  for any rigid motion impressed initially to the coupled system $\mathscr S$,  the asymptotic motion will be a permanent rotation around a central axis of inertia,   aligned with the total initial angular momentum $\bfK_G$. Moreover, taking into account that, in this case,
$$
A=B=\Frac{M}{12}(3R^2+h^2)\,,\ \ C=\Frac{1}{2}MR^2\,,
$$
with $M$ mass of the liquid, 
the axis of rotation coincides with ${\sf a}$ if $h<\sqrt{3}R$, whereas if $h>\sqrt{3}R$ it is  perpendicular to ${\sf a}$, and passing through $G$. Therefore, the rotation occurs along the axis of the cylinder if and only if the cylinder is ``sufficiently flattened''; see Figure 1.   \bigskip\par
\ER{4}
\begin{figure}
\vspace*{0.009in}
\begin{center}
\def\svgwidth{\textwidth}
\scriptsize{
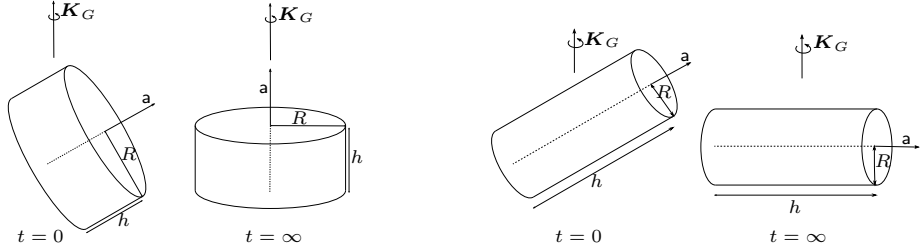}
\end{center}
\caption{
Initial and final orientation of a cylinder filled with a viscous liquid when $h<\sqrt{3}R$ (left), and   $h>\sqrt{3}R$ (right).  }
\end{figure}
\Br Results proved in \theoref{4} require the initial data to be in a certain range (see \Eqref{8.3}, \Eqref{8.5}, and \Eqref{PN2}). However, the numerical tests reported in Section 9.2, suggest that such a requirement might be unnecessary. The question of whether this restriction can be removed analytically is at the moment open. 
\ER{4+}
\smallskip\par
With the help of \theoref{4}  we are now able to derive the following results, which ensure    stability  of permanent rotations of the coupled system around the central axis with the largest moment of inertia, and instability in the other cases. 

\Bt Let  $\mathscr S$ perform a permanent rotation around the  central axis $\{G,\bfe\}$, namely, $\bfv\equiv \0$,    $\bfomega_\infty\equiv\bfomega_\infty^{(0)}=\omega_0\bfe$, $\bfe\in\{\bfe_1,\bfe_2,\bfe_3\}$. Assume that, at time $t=0$, this state is perturbed, and denote by $\bfv=\bfv(x,t)$, $\tilde{\bfomega}_\infty(t)=(\tilde p(t),\tilde q(t),\tilde r(t))$ the corresponding perturbation fields.

The following properties hold.
\begin{itemize}
  \item [{\rm (a)}] If $A<B=C$, then the  permanent rotation with $\bfe\equiv\bfe_1$ is unstable in the sense of Lyapunov.
  \item[{\rm (b)}] If $A\le B<C$,  then the  permanent rotation with $\bfe$ being either $\bfe_1$ or $\bfe_2$ is unstable in the sense of Lyapunov. If, however, $\bfe\equiv\bfe_3$, then the corresponding permanent rotation is   stable. Precisely,
for any $\varepsilon>0$ there is $\delta>0$ such that
\be
E(0)+\tilde p^2(0)+\tilde q^2(0)+\tilde r^2(0)<\delta\ \ \Longrightarrow\ \ E(t)+\tilde p^2(t)+\tilde q^2(t)+\tilde r^2(t)<\varepsilon\,,
\eeq{8.8} 
\mbox{for all $t\ge 0$}. Moreover, there is $\gamma=\gamma(A,B,C,\omega_0)>0$ such that if
\be
E(0)+ \tilde p^2(0)+\tilde q^2(0))+ \tilde r^2(0)\le \gamma\,,
\eeq{8.9}
it results\footnote{Recall that, by \theoref{3}, $E(t)\to 0$, as $t\to\infty$ regardless of the ``size'' of the initial conditions.} 
 \be \tilde p(t),\tilde q(t)\to 0\,, \ \ \tilde r(t)\to { r^*},\ \ \mbox{as $t\to\infty$}\,,
\eeq{gmb}
where
\be
 r^*=-\omega_0\pm\sqrt{\frac{1}{C^2}(A^2\tilde p^2(0)+B^2\tilde q^2(0)^2)+(\tilde r(0)+\omega_0)^2}\,,
\eeq{8.10}
and where we take $+$ or $-$ according to whether $\omega_0>0$ or $\omega_0<0$.
\item[{\rm (c)}] If $A=B=C$, the permanent rotation corresponding to arbitrary $\bfe$ is stable in the sense of Lyapunov,  namely, \Eqref{8.8} holds
\mbox{for all $t\ge 0$}.
\end{itemize}
\ET{5}
{\bf Proof.}  We begin to notice   that $(\bfv,\tilde\bfomega_\infty+\bfomega_\infty^{(0)})$   must satisfy \Eqref{8.1} and \Eqref{8.2}, and, consequently, we may apply \theoref{4}. In order to show the property in (a), take initial conditions for the perturbed field  satisfying $E(0)=\tilde q(0)=\tilde r(0)=0$, and $\tilde p(0)$ non-zero and as small as we please. Thus, in particular, \Eqref{8.5} is satisfied. As a consequence, by (b) of \theoref{4} and \Eqref{PN5}$_2$, we must have
$$
B^2(\bar q^2+\bar r^2)= A^2(\omega_0^2+\tilde p^2(0))\,.
$$
Thus, for all $\tilde p(0)$ sufficiently small, 
there is $\bar t>0$ independent of $\tilde p(0)$, such that $(\tilde q^2(t)+\tilde r^2(t))\ge \half (A^2 \omega_0^2/B^2)$, for all $t\ge \bar t$, which furnishes the desired instability result. Next, to prove the first property stated in (b),  we take (in both cases $\bfe=\bfe_1,\bfe_2$) $E(0)= \tilde p(0)=\tilde q(0)=0$ and $\tilde r(0)$ arbitrarily small, and notice that \Eqref{8.3} and \Eqref{PN2} is satisfied. By a completely analogous reasoning to the one employed previously we then show $\tilde r^2(t) \ge \half (A^2 \omega_0^2/C^2)$ for all sufficiently large $t$, thus proving instability. To show the other statement in (b), we multiply both sides of \Eqref{8.1} by $C$, and subtract to the resulting inequality \Eqref{8.2}, side by side. We deduce, in particular,
$$\ba{rl}\medskip
2CE(t)+A(C-A)\tilde p^2(t)&\!+B(C-B)\tilde q^2(t)\\
&\le 2CE(0)+A(C-A)\tilde p^2(0)+B(C-B)\tilde q^2(0)\,,
\ea
$$ 
which, in turn, implies
\be\ba{cc}\medskip 
E(t)+\tilde p^2(t)+\tilde q^2(t)\le m\big[E(0)+\tilde p^2(0)+\tilde q^2(0)\big]\,,\\ m:=\Frac{\max\{2C,A(C-A),B(C-B)\}}{\min\{2C,A(C-A),B(C-B)\}}\,.\ea
\eeq{mg1}
Thus,  given $\varepsilon>0$, we have
\be
E(0)+\tilde p^2(0)+\tilde q^2(0)<\Frac{\delta_1}{m}\ \Longrightarrow\ \ E(t)+\tilde p^2(t)+\tilde q^2(t)< {\delta_1}\,,\ \ \mbox{for all $\delta_1\in (0,\varepsilon/2)$\,.} 
\eeq{mg2}
Next, we want to show that for a suitable choice of $\delta_2>0$, the following property holds
\be
\tilde r^2(0)<\delta_2\ \ \Longrightarrow\ \ \tilde r^2(t)<\frac{\varepsilon}{2}\,.
\eeq{mg3}
Without loss of generality, we take \be\varepsilon=2\eta^2\omega_0^2\,,  \eeq{mg4}\mbox{with $\eta$  arbitrarily fixed in $(0,1)$}, and choose $\delta_2<\varepsilon/2$. Assume \Eqref{mg3} is not true. In view of the continuity of $r(t)$, let $\bar t>0$ be the {\em first} instant of time such that $\tilde r^2(\bar t)=\varepsilon/2$. Thus, by \Eqref{8.2} and \Eqref{mg4} we deduce
\be 
\pm C^2\omega_0^2\eta(2\pm\eta)=A^2\tilde p^2(0)+B^2\tilde q^2(0)+C^2\tilde r(0)(\tilde r(0)+2\omega_0)-A^2p^2(\bar t)-B^2\tilde q^2(\bar t)\,.
\eeq{mg5}
Recalling that $A\le B$, $\eta\in (0,1)$ and using \Eqref{mg2}, from \Eqref{mg5} we show
$$
C^2\omega_0^2\eta\le B^2\left(\Frac{m+1}{m}\right)\delta_1+C^2\sqrt{\delta_2}\big(\sqrt{\delta_2}+2|\omega_0|\big)\,.
$$
Employing in the latter relation the inequality
$2\sqrt{\delta_2}|\omega_0|\le 2\delta_2/\eta+\eta\omega_0^2/2$, and recalling again that $\eta\in (0,1)$, we get
\be
C^2\omega_0^2\eta \le 2B^2\left(\frac{m+1}m\right)\delta_1+6C^2\frac{\delta_2}{\eta}\,.
\eeq{mg6}
However, \Eqref{mg6} cannot be true as long as we pick $\delta_1,\delta_2$ such that (for instance)
$$\ba{ll}\medskip
0<\delta_1<\Frac{m\,C^2\omega_0^2}{2(m+1)B^2}\Frac{\eta}{4}\equiv \Frac{m\,C^2|\omega_0|}{8(m+1)B^2}\Frac{\sqrt{\varepsilon}}{\sqrt{2}} 
\\
0<\delta_2<\omega_0^2\Frac{\eta^2}{24}\equiv \Frac{\varepsilon}{48} \,.\ea
$$ 
As a consequence, \Eqref{8.8} follows, provided we choose
$$
\delta<\min\Big\{\Frac{\varepsilon}{48}\,,\Frac{\varepsilon}{2m}\,, \Frac{m\,C^2|\omega_0|}{8(m+1)B^2}\Frac{\sqrt{\varepsilon}}{\sqrt{2}}\Big\}\,. 
$$
Let us now show  the last property stated  in (b). From \theoref{4}  we know that the  asymptotic property \Eqref{gmb}  is valid  
whenever the initial conditions of the motion $(\bfv,\tilde\bfomega_\infty+\bfomega_\infty^{(0)})$ satisfy \Eqref{8.3} and \Eqref{PN2}. Recalling that $A\le B<C$,  one shows that both conditions are certainly met if
$$
E(0)+\frac{A}{2B}(B-A)(\tilde q^2(0)+\tilde p^2(0))\le \frac{C}{2B}(C-B)(\omega_0+\tilde r(0))^2\,.
$$
However, since $(\tilde r(0)+\omega_0)^2\ge \half \omega_0^2-\tilde r^2(0)$, we see that the latter is satisfied provided $E(0)$ and $\tilde p(0),\tilde q(0)$, and $\tilde r(0)$ obey \Eqref{8.9}, for a suitable definition of $\gamma$. However, by taking $\gamma$ even smaller if necessary, from the stability property  proved above we know that $|\tilde r(t)|<|\omega_0|$, for all $t\ge 0$. Consequently, \Eqref{8.10} follows from this consideration, by passing to the limit $t\to\infty$ in \Eqref{8.2}.  It remains to show property (c). In this regard, we observe that from our hypothesis and \Eqref{3.8}$_3$ we deduce
\be
\dot{\bfomega}_\infty=\bfa\times(\bfomega_\infty+\bfomega_\infty^{(0)})\,,
\eeq{SF3}
from which it follows that
\be
|\bfomega_\infty(t)+\bfomega_\infty^{(0)}|=
|\bfomega_\infty(0)+\bfomega_\infty^{(0)}|\,.
\eeq{SF4}
From \Eqref{SF3} and \Eqref{SF4} we thus obtain
\be
|\bfomega_\infty(t)|\le |\bfomega_\infty(0)|+|\bfomega_\infty(0)+\bfomega_\infty^{(0)}|\int_0^\infty|\bfa(t)|\,.
\eeq{SF5}
Using Schwartz inequality and \Eqref{SFX} in \Eqref{SF5} allows us to conclude
$$
|\bfomega_\infty(t)|\le |\bfomega_\infty(0)|+c\,|\bfomega_\infty(0)+\bfomega_\infty^{(0)}|\,\|\bfv(0)\|_2\,,
$$
so that the property stated in (c) follows from this last inequality and \Eqref{SFX}.
\QED
\begin{remark} Combining  \theoref{4} and \theoref{5} we derive the following interesting consequence. Suppose  $C>A\ge B$, $E(0)=0$, $p(0),q(0)$ sufficiently ``small'', and $r(0)\neq 0$. Then,  the asymptotic motion of the coupled system  --which we know is a permanent rotation around $\bfe_3$-- will have   angular velocity $\bar\bfomega=\bar r\bfe_3$, where $\bar r$ has the {\em same sign as $r(0)$}. Observing that  $\bar\bfomega= \kappa\,\bfK_G$, $\kappa>0$,  this property implies that $\{G,\bfe_3\}$ has to keep (asymptotically) the same {\em orientation} with  $\bfK_G$ that it had at time $t=0$; see Figure 2.  Stated differently, this means that, at least under the above conditions, the axis $\bfe_3$ cannot (eventually) ``flip-over''. This property is confirmed by the numerical tests presented in Section 9.2; see Figure 5, bottom panel. These tests also show, however,   that the above property is no longer valid for initial data of {\em finite} size; see Figure 5, top panel. Similar experiments  prove, in addition,  that if $r(0)=0$, a change of viscosity of the liquid may  trigger such an effect as well. In other words, it is found that in some range of viscosities the orientation of $\bfe_3$ and $\bfK_G$ is the same, whereas in another range it is opposite; see Section 9.3. 
It will be the object of  future work to investigate the analytical aspect of this interesting phenomenon.
\label{rem:6}
\end{remark}

\begin{remark}The instability result in \theoref{5} should be contrasted  with their ``classical'' counterpart when the cavity in the gyroscope is empty. In fact, in such a situation, as is well known, permanent rotations  are stable in {\em both} cases $A,B\,\substack{ < \\ > }\, C$; see \cite[Chapter 8.33]{Lei}. Whereas, if the cavity is filled with a viscous liquid, this permanent rotation is (axially) stable if and only if $C>A,B$. 
\end{remark}
\bigskip\par
\begin{figure}
\begin{center}
\def\svgwidth{
2in}
\scriptsize{
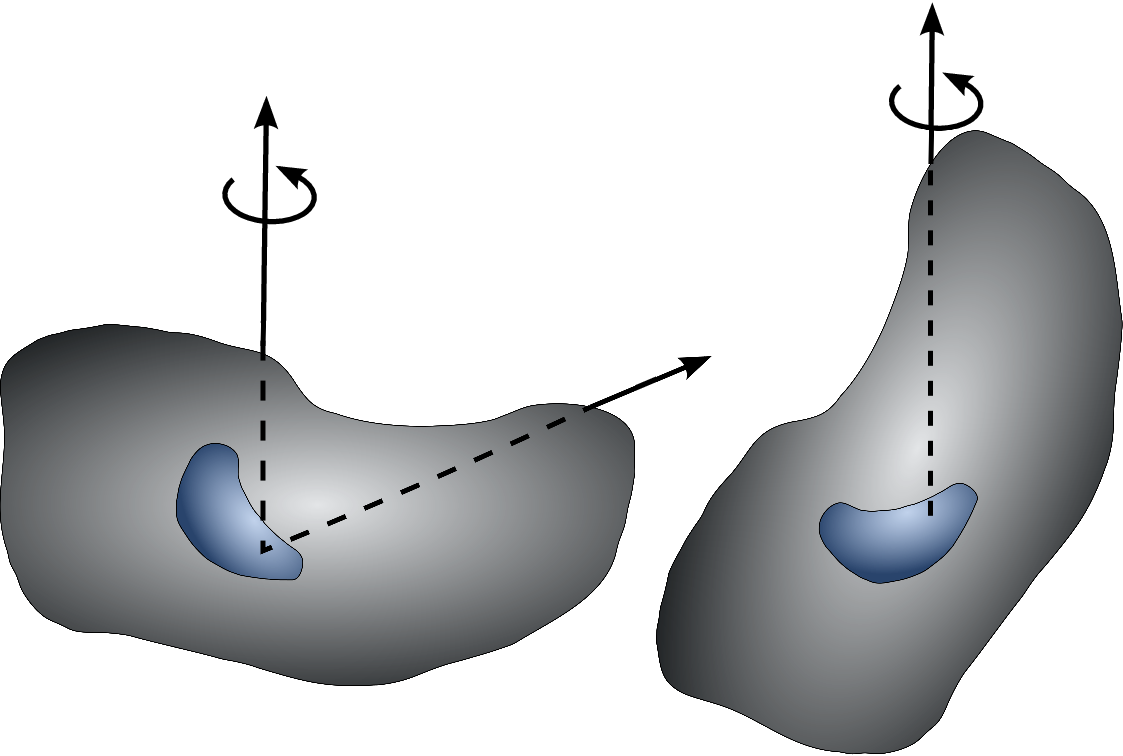}
\def\svgwidth{
2in}
\hspace*{0.38in}
\scriptsize{
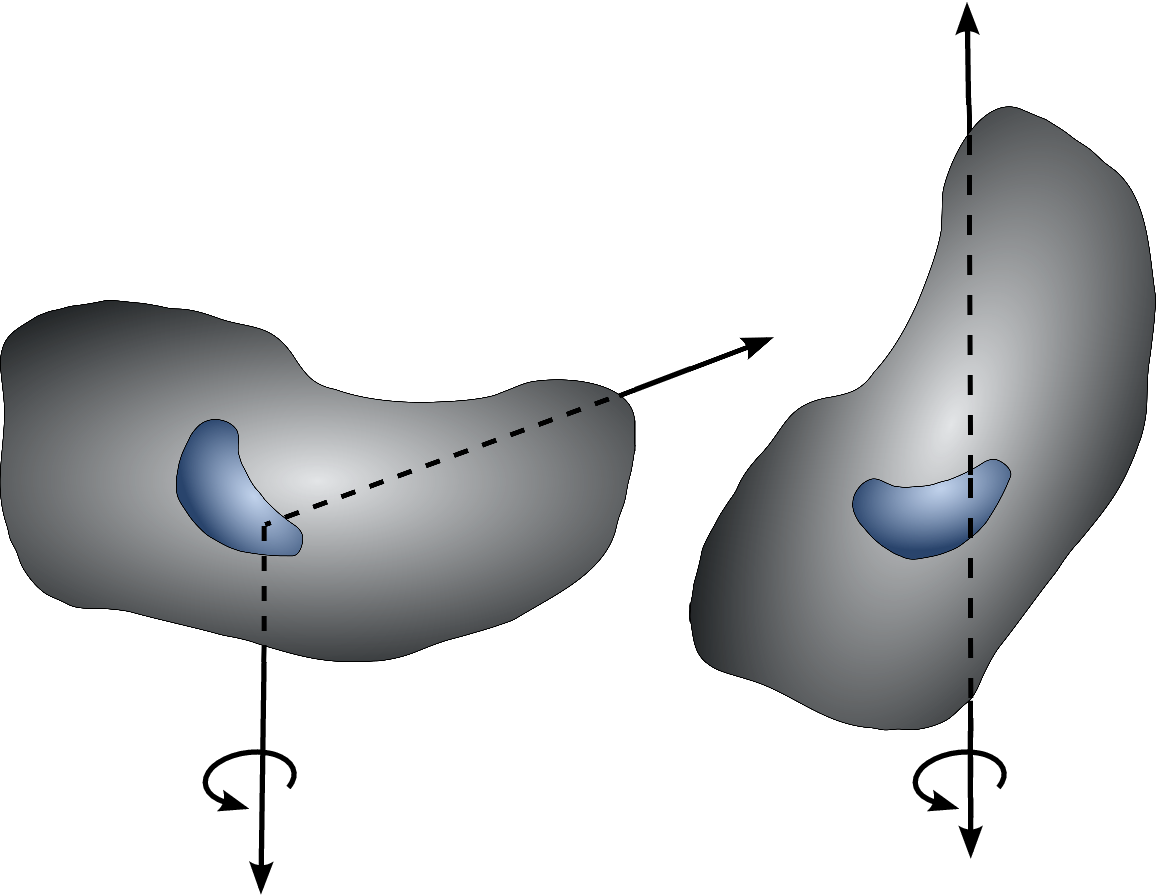}
\end{center}
\caption{
Dependence of the orientation of a body filled with a viscous liquid on the initial angle $\theta(0)$ between $\bfe_3$ and $\bfK_G$: $\theta(0)<\pi/2$ (left);   \ $\theta(0)>\pi/2$ (right). }
\end{figure}

\def\body{{\mathscr{B}}}
\def\fluid{{\mathscr{C}}}
\def\SS{\mathcal{S}_\body}
\def\yy{{\boldsymbol y}}
\def\xx{{\boldsymbol x}}
\def\nn{{\boldsymbol n}}
\def\inertiab{\bfI}
\def\inertiabf{\inertiab}
\def\vel{{\boldsymbol u}}
\def\velb{\boldsymbol{\xi}}
\def\relv{{\boldsymbol v}}
\def\relvarg{\relv(\vel,\angvel,\velb)}
\def\angvel{\boldsymbol{\omega}}
\def\cauchy{\mathbf{T}}
\def\pres{{\sf p}}
\def\cauchyarg{\mathbf{T}(\vel,\pres)}
\def\FF{\mathbf{F}}
\def\Pol{\mathbb{P}}
\section{Numerical experiments}

In order to complement and, to an extent,  complete the analytical results described in the previous sections, we have performed a number of targeted numerical simulations, also as a tool to provide insight and further information on the long term behavior of \Eqref{3.1}. As we know, conservation of angular momentum and  energy balance are the fundamental governing principles for the system at hand. It is, therefore, essential that the time discretization method accurately preserves the invariants of the system. At the level of numerical approximation this is not a trivial task, even when we restrict ourselves to the analysis of the motion of the body solely, i.e. equations \Eqref{3.1}$_3$. The dynamics of a rigid body can be formulated in the framework of Hamiltonian mechanics. The family of \emph{geometric} or \emph{symplectic integrators} are particularly adapted to approximate such problems, because they preserve the invariants of the Hamiltonian flow \cite{MR2221614}. On the basis of these considerations, we have employed the following time discretization algorithm for the coupled liquid / body system whose motion is governed by \Eqref{3.1}. For the time integration of the body dynamics, we use the $\theta$-method. In particular, we adopt the implicit midpoint integration rule, $\theta=\frac12$, because of its good properties as a geometric integrator. More precisely, it is known that the midpoint rule exactly satisfies the conservation of momentum for simple Hamiltonian systems. For the Navier-Stokes equations, we apply the implicit Euler time advancing scheme in order to guarantee the stability of the algorithm. At each time step, we use sub-iterations to uncouple the solution of the discrete body and liquid problems and to linearize the corresponding equations. The combination of these techniques gives rise to the following algorithm:

\begin{description}
\item given $\vel_0,\,\angvel_0$ and a partition of the interval $(0,T]$ in evenly distributed time steps $t_n=n \tau$ with $\tau >0$, for $n=1,2,3,\ldots$ find $\vel_n,\,\pres_n,\,\angvel_n$ such that:

\medskip

\item given $\vel_n^0=\vel_{n-1},\,\angvel_n^0=\angvel_{n-1}$, for $k=1,2,3,\ldots$ solve the following sub-problems,

\medskip

\item body problem: find $\angvel_n^*$ such that,
\begin{multline*}
\inertiab \tau^{-1} \big(\angvel_n^*-\angvel_{n-1}\big) + \theta \big[ \angvel_n^{k-1} \times \big(\inertiab \cdot \angvel_n^{k-1} \big) + \int_{\partial\fluid} \xx \times \bfT (\vel_n^{k-1},\pres_n^{k-1})\cdot \nn \big] 
\\
+(1-\theta) \big[ \angvel_{n-1} \times \big(\inertiab \cdot \angvel_{n-1} \big) + \int_{\partial\fluid} \xx \times \bfT (\vel_{n-1},\pres_{n-1})\cdot \nn  \big] = \0
\end{multline*}

\item relaxation: given $\sigma \in (0,1)$ calculate $\angvel_n^k = \sigma \angvel_n^* + (1-\sigma) \angvel_{n-1}$

\medskip

\item liquid problem: find $\vel_n^k,\,\pres_n^k$ such that,
\begin{equation*}
\begin{cases}
\rho  \big(\vel_n^k - \vel_{n-1} \big) + \rho \angvel_n^k \times \vel_n^k + \rho\,  \relv(\vel_n^{k-1},\angvel_n^k) \cdot \nabla  \vel_n^k 
\\
\hfill - \Div \bfT(\vel_n^k,\pres_n^k) = \0 \quad  \text{in} \ \fluid, 
\\
\Div \vel_n^k = 0, \quad  \text{in} \ \fluid,
\\
\vel_n^k = \angvel_n^k \times \xx, \quad \text{on} \ \partial\fluid,
\end{cases}
\end{equation*}

\item convergence test: 
given $\epsilon$ small enough, if $\|\angvel_n^k - \angvel_n^{k-1}\| < \epsilon$ then $\vel_n=\vel_n^k,\,\pres_n=\pres_n^k,\,\angvel_n=\angvel_n^k$.
\end{description}

\begin{figure}
\begin{center}
\includegraphics[width=0.5\textwidth]{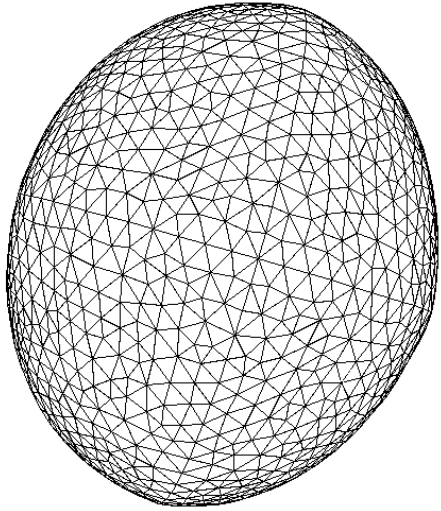}
\caption{Geometrical configuration of the cavity  $\fluid$ used for the numerical experiments.}
\label{fig:geometry}
\end{center}
\end{figure}

For the spatial approximation of the liquid equations we exploit the finite element method \cite{MR851383,MR1017032,MR2454024}. We address the saddle point formulation of the problem in terms of velocity and pressure variables. In order to achieve a stable discretization of the divergence-free constraint, we use \emph{inf-sup} stable mixed finite elements, such as $\Pol^2-\Pol^1$ approximation of the velocity and pressure fields, respectively, see \cite{MR1115205,MR851383}. Since our numerical tests involve relatively simple geometrical configurations, moderately refined computational grids will be applied, see for instance Figure \ref{fig:geometry} that shows the geometrical model of $\fluid$ and the corresponding computational mesh. The system of algebraic equations arising from the discretization scheme is solved  by means of direct techniques, which turn out to be an effective option since the number of degrees of freedom is not excessively large.

\subsection{Effect of the  viscosity on the final angular velocity}
\begin{figure}
\begin{center}
$\nu=0.1,\, \angvel_0 = [6.2697,0.4109,0],\,E(0)=0,\, t\in(0,80)\, s$\\ \includegraphics[width=\textwidth]{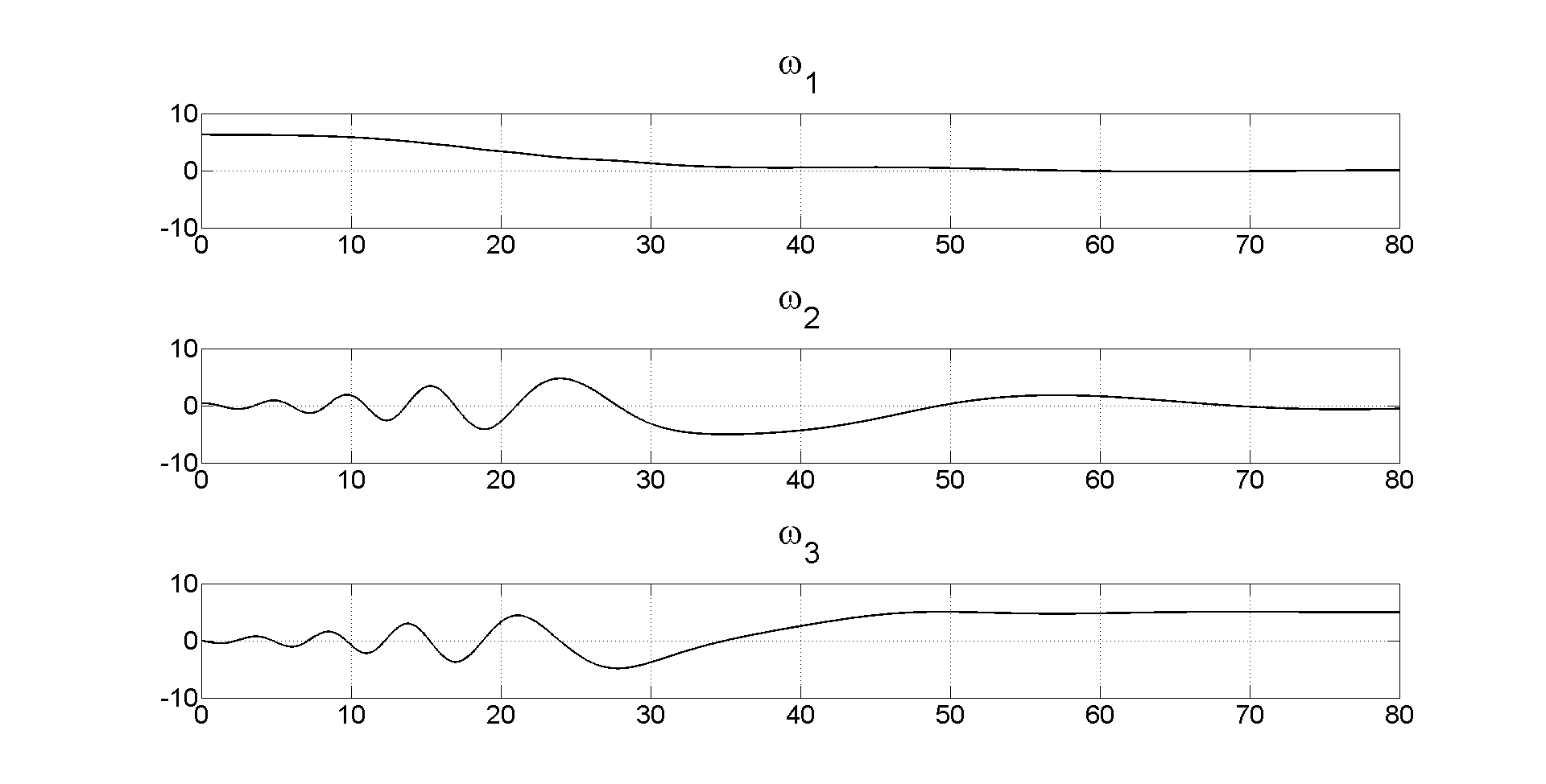}\\
$\nu=0.001,\, \angvel_0 = [6.2697,0.4109,0],\,E(0)=0,\, t\in(0,80)\, s$\\ \includegraphics[width=\textwidth]{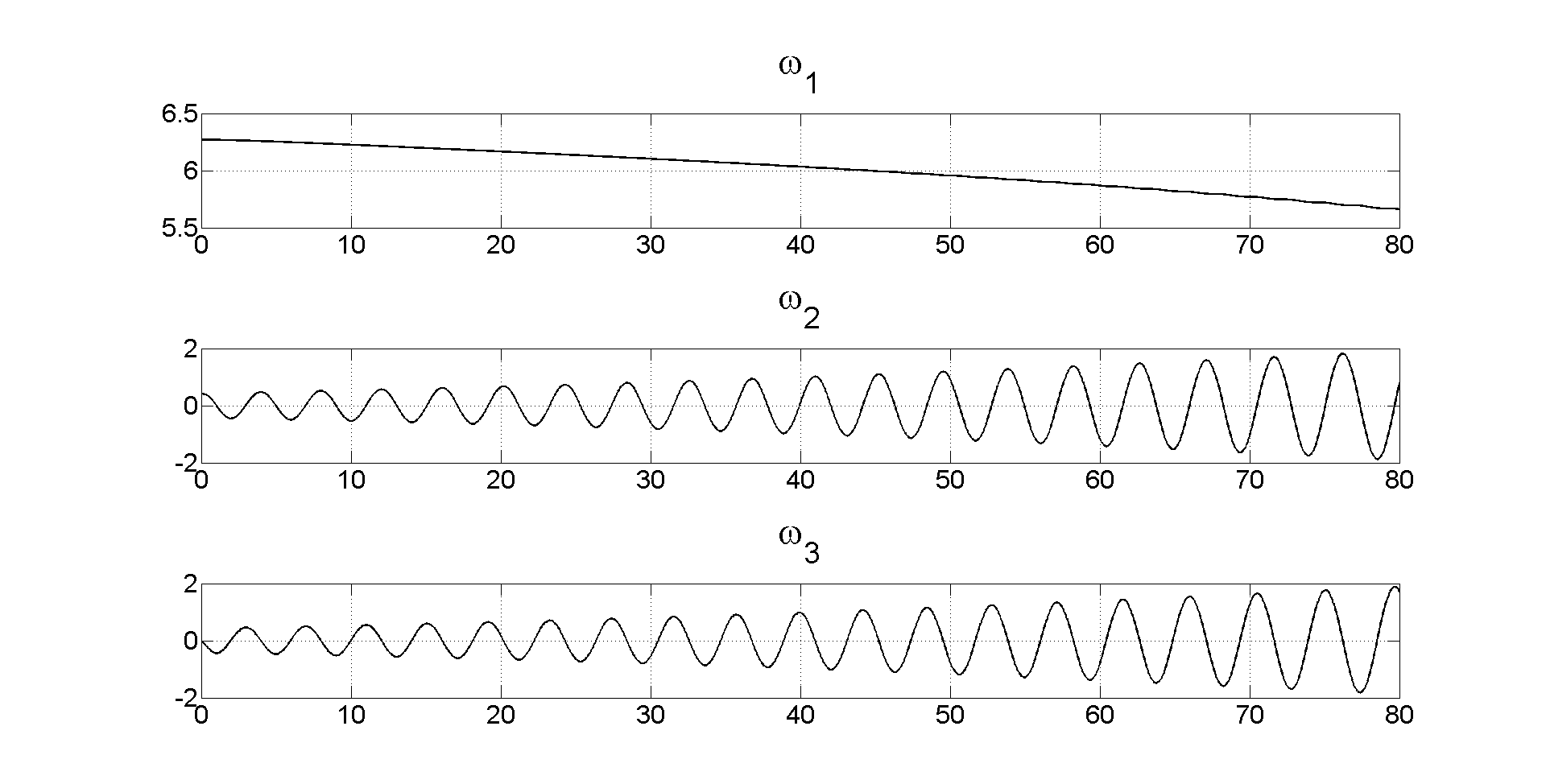}\\
$\nu=0.001,\, \angvel_0 = [6.2697,0.4109,0],\,E(0)=0,\, t\in(0,1200)\, s$\\ \includegraphics[width=\textwidth]{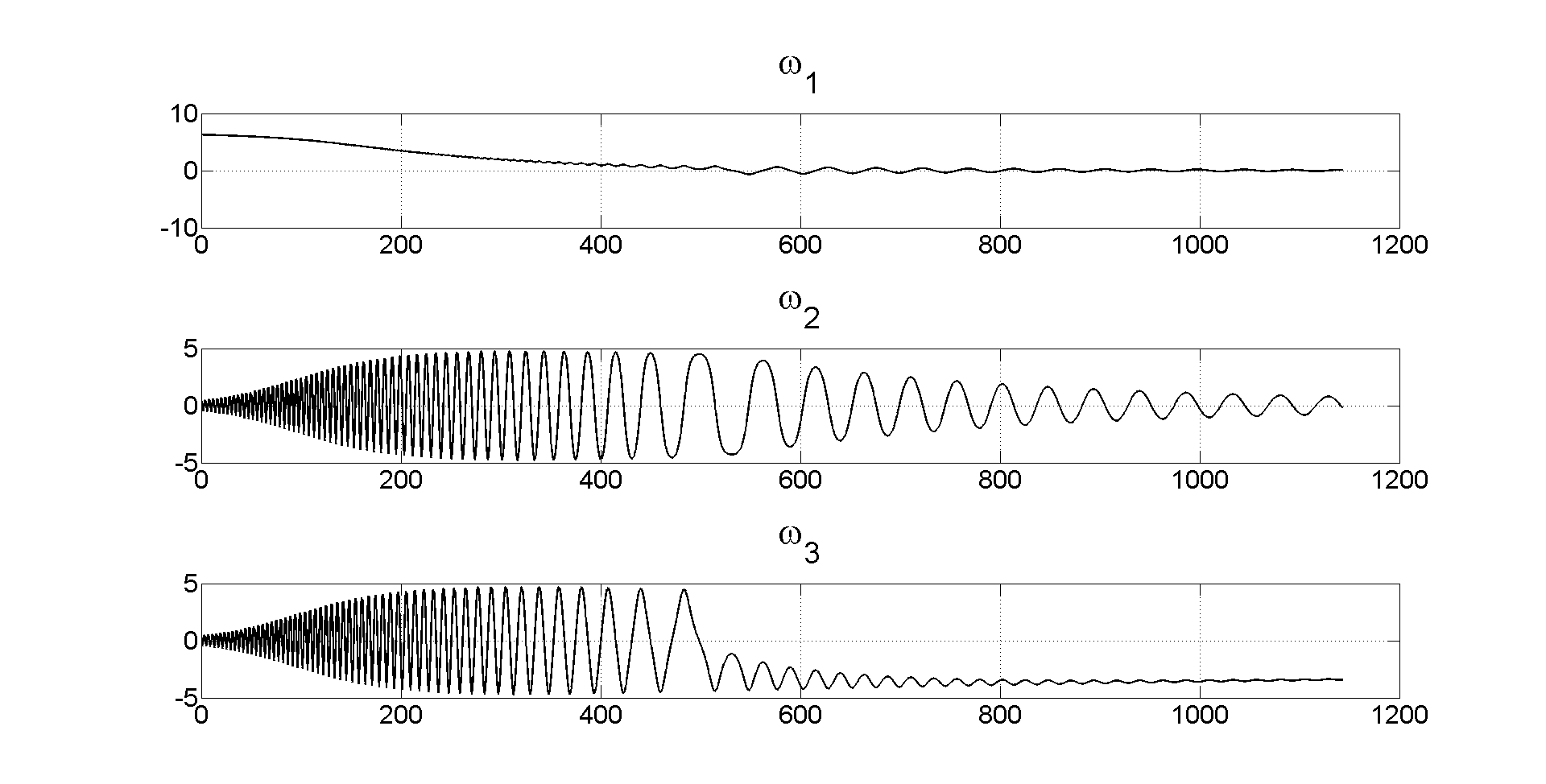}
\caption{Dynamics of the liquid-solid system for decreasing values of the viscosity.}
\label{fig:test0}
\end{center}
\end{figure}

We study the dynamics of a system where the cavity has the (quasi-ellipsoidal) shape shown in Figure \ref{fig:geometry} and the rigid body is characterized by a tensor of inertia with eigenvalues $A=5.54,\, B=6.73,\, C=6.76$ (which corresponds to depositing a layer of uniform material of constant thickness around the cavity). At the initial time, the motion of $\body$ is identified by the angular velocity $\angvel(0)=2\pi[\cos(\theta), \cos(\phi)\sin(\theta), \sin(\phi)\sin(\theta)]$ with $\theta=\pi/48,\,\phi=0$, while the relative velocity of the liquid is $\boldsymbol v=\boldsymbol 0$ everywhere in $\fluid$. In Figure \ref{fig:test0} we visualize the plots of $\angvel(t)=[p(t),q(t),r(t)]$ for decreasing values of the kinematic viscosity of the liquid, namely $\nu=\mu/\rho$ as in equation \Eqref{3.1}. The numerical simulations show that, for moderately large values of the viscosity ($\nu=0.1$), the system quickly reaches a steady rotation around  the central axis of inertia with the larger moment of inertia. As expected, the trend through which the rotational equilibrium is reached is extremely sensitive to the magnitude of the viscosity. Indeed, for $\nu=0.001$ the rotation of the liquid-solid system is ``chaotic", at least for the timescale of analysis appropriate for large viscosity. Only when timespan of simulation is significantly extended, the numerical experiments show that the steady rotation is eventually recovered (Figure \ref{fig:test0}, bottom panel).
The numerical simulations enable a more quantitative analysis of the effect of the viscosity on the time required to reach equilibrium. Let us denote by $t_c$ the instant at which the following condition is satisfied for the first time,
\[t_c : \quad \frac{\|\angvel(t) - \angvel(\infty)\|}{\|\angvel(0) - \angvel(\infty)\|} < 0.1\]
For the viscosities $\nu=0.1,\,0.05,\,0.02,\,0.01$ we have calculated $t_c$ and the corresponding $\angvel(t_c)$. The results are reported in Table \ref{tab:1}. From these data, it is possible to estimate how $t_c$ depends on $\nu$. We begin by postulating a power law dependence such as $t_c \simeq \nu^p$. Then, in the range $\nu \in (0.01,0.1)$ the value of $p$ that best fits the data is $p=-0.305$. This result confirms the inverse dependence of the time required to reach equilibrium on the magnitude of the viscosity.
\begin{table}
\begin{center}
\begin{tabular}{c|ccc|c}
$\nu$ & $p(t_c)$ & $q(t_c)$ & $r(t_c)$ & $t_c$\\\hline
0.1	&	-0.4018	&	0.4558	&	4.8682	&	50.8\\
0.05	&	0.4908	&	-0.568	&	4.7584	&	63.3\\
0.02	&	-0.4887	&	-0.4171	&	-4.5768	&	75.2\\
0.01	&	0.6319	&	0.4658	&	-4.3192	&	99.8\\
\end{tabular}
\end{center}
\caption{Dependence of the (numarically estimated) time to reach equilibrium on the liquid kinematic viscosity. The initial rotation is $\angvel_0 = [6.2697,0.4109,0]$.}
\label{tab:1}
\end{table}

\subsection{Effect of the initial rotation on the final angular velocity}
\begin{figure}
\begin{center}
%
$\nu=0.1,\, \angvel(0)=[4.44,3.14,3.14],\, E(0)=0$\\
\includegraphics[width=\textwidth]{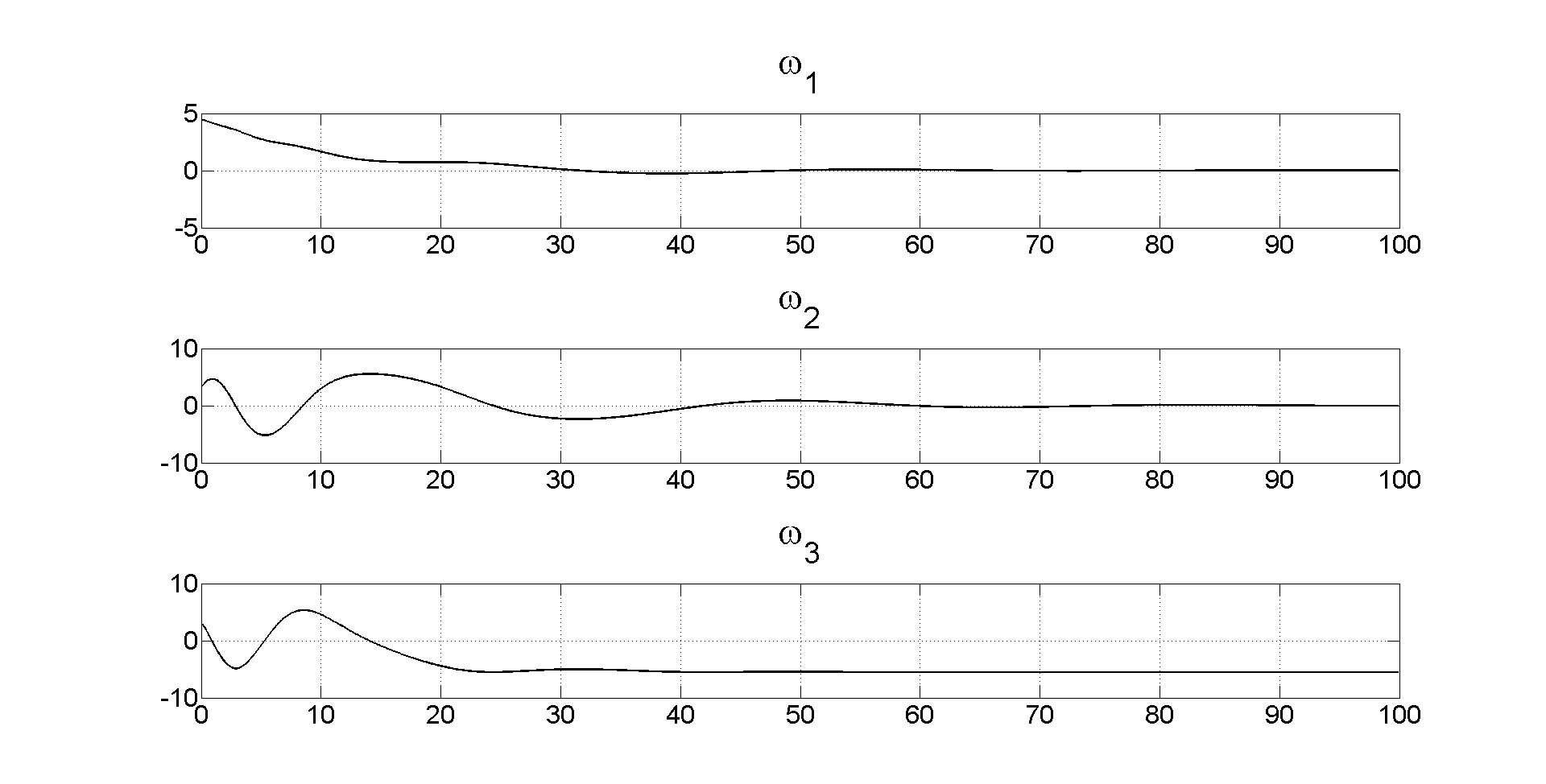}\\
%
$\nu=0.1,\, \angvel(0)=[0.444,0.314,3.14],\, E(0)=0$\\
\includegraphics[width=\textwidth]{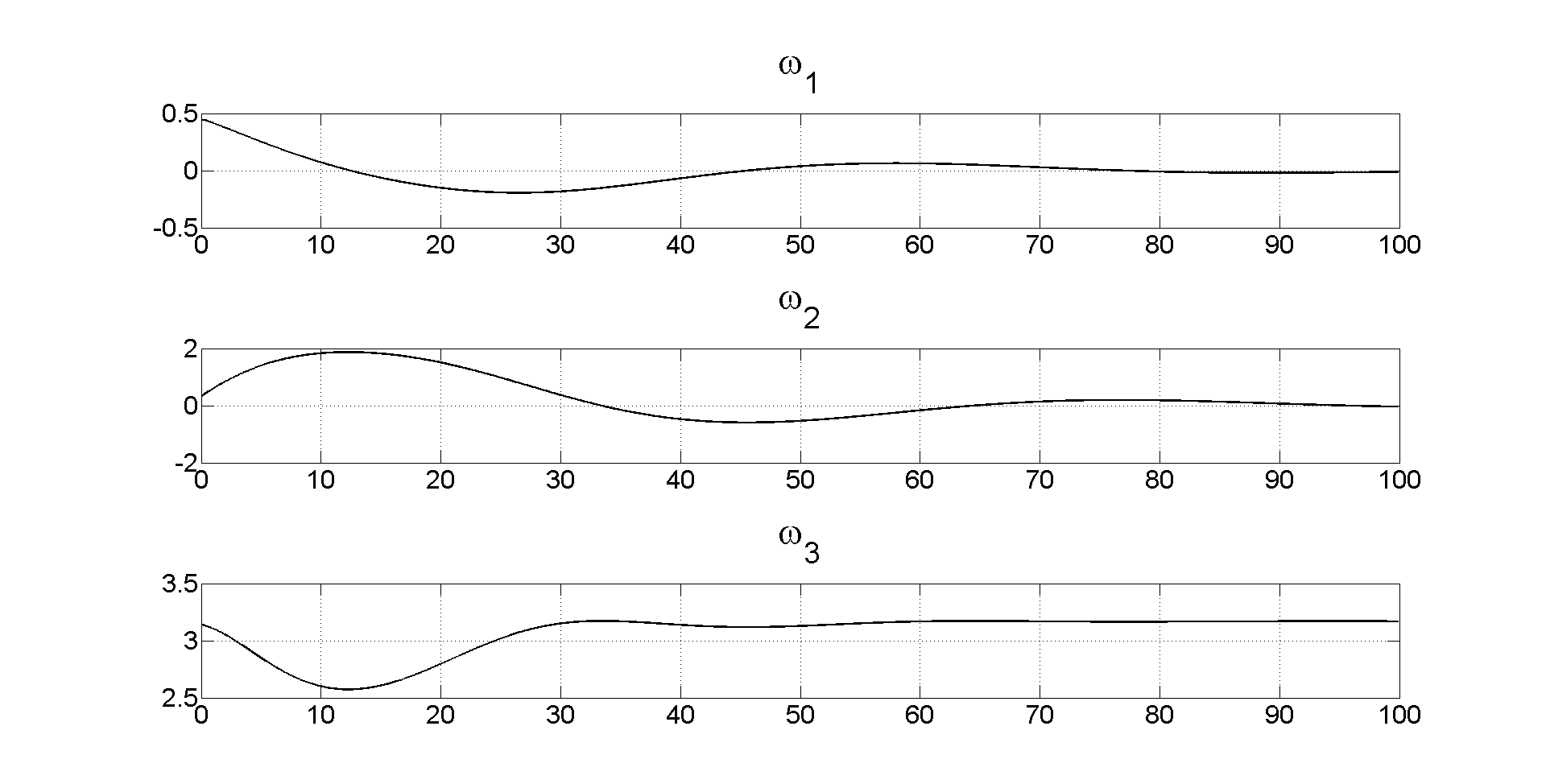}
\caption{Dynamics of the liquid-solid system for different inital angular velocity.}
\label{fig:test1}
\end{center}
\end{figure}

\begin{figure}
\begin{center}
$\nu=0.1,\, \angvel(0)=[4.44,3.14,3.14],\, E(0)=0$\\
\includegraphics[width=\textwidth]{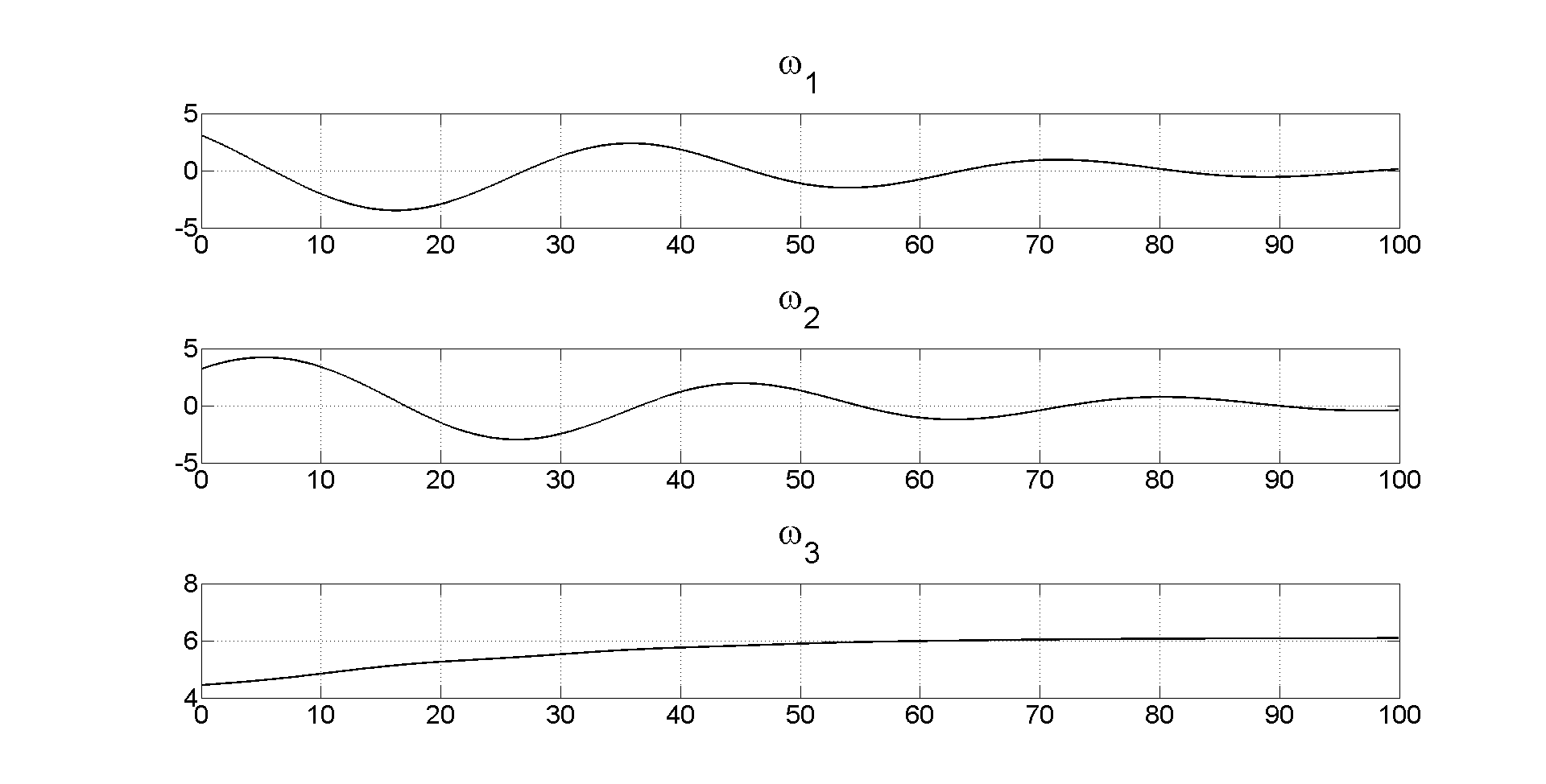}\\

$\nu=0.1,\, \angvel(0)=[4.44,3.14,3.14],\, E(0)=341.6$\\
\includegraphics[width=\textwidth]{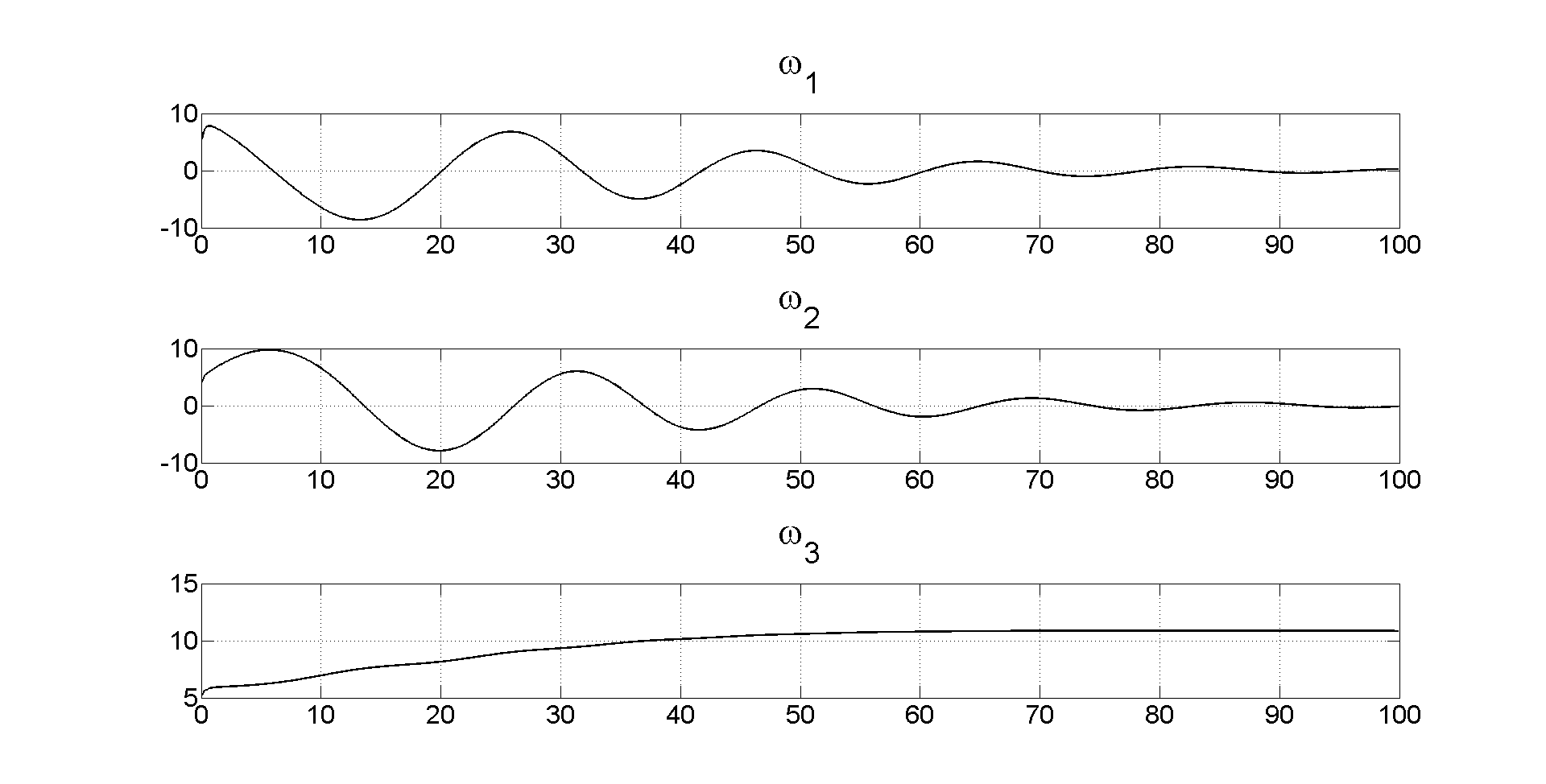}\\
\caption{Dynamics of the liquid-solid system for different initial liquid energy $E(0)$.}
\label{fig:test3}
\end{center}
\end{figure}

We discuss here experiments to investigate the attainability of permanent rotations, for which the main analytical result are Theorem 4 and Theorem 5. 
The numerical experiments are particularly helpful to test the validity of the analysis beyond the restrictions \Eqref{8.3}, \Eqref{8.5}, \Eqref{PN2}. 
In these cases, the liquid kinematic viscosity is set to $\nu=0.1$. 
We begin with a problem configuration where the condition \Eqref{PN2}(a) is not satisfied because of large initial rotational speed,
\[
9.6860=\frac{A}{2B}(B-A)p^2 > \frac{C}{2B}(C-B)r^2 = 0.1310.
\]
The computed plots of $\angvel(t)$ reported in Figure \ref{fig:test1} (top panel) show that,  
for ``large" initial data, the conclusions of Remark 9 are no longer valid. In particular, we observe that in this case $r(0)>0>\overline{r}$, which differs from the predicted behavior for small data, while $\overline{p}=\overline{q}=0$ as proved in the analysis. In other words, sufficiently large $p(0), q(0)$ may trigger the ``flip-over" effect. In Figure \ref{fig:test1} (bottom panel) we investigate a similar situation, where $p(0)$ and $q(0)$ are sufficiently small. In this case, the validity of condition \Eqref{PN2}$_1$ is restored and the results of Remark 9 (i.e. $r(0)$ and $\overline{r}$ share the same sign) are reproduced, as expected, by the numerical simulation.
Numerical experiments also elucidate the behavior of the system when the initial relative velocity of the liquid with respect to the rigid body is varied. More precisely, we compare two cases which only differ from the initial liquid energy $E(0)$. In one case the initial   relative velocity of the liquid is initialized to $\boldsymbol v = \boldsymbol 0$ in $\fluid$. As a result $E(0)=0$. In the other case  we define  $\boldsymbol v$ as a nonzero compatible velocity field, such that $\boldsymbol v = \boldsymbol 0$ on $\partial\fluid$ and $\nabla \cdot \boldsymbol v =0$ in $\fluid$, such that $E(0) \gg 0$. In particular, at the initial time the (absolute) liquid velocity in the reference frame attached to $\body$ can be expressed in the following form $\boldsymbol u = f(\|\boldsymbol x\|) \angvel(0) \times \boldsymbol x$. Since $f(\|\boldsymbol x\|) \neq 1$ then $\boldsymbol v \neq 0$. For this numerical experiment we consider a different tensor of inertia, $A=4.99,\, B=4.99,\, C=5.54$. In the case with $E(0) \gg 0$ we observe that \Eqref{8.3} is not satisfied because
\[ 6.0786 = \frac{(C-A)C}{2A}r^2(0) \ll E(0) = 341.6 \]
The corresponding results are shown in Figure \ref{fig:test3}. We see that the qualitative behavior of the system is substantially unaffected. This test shows that the thesis of Theorem 4, case (a) may be verified also when condition \Eqref{8.3} is violated.  This fact may be interpreted observing that for large viscosities, the liquid quickly adjusts its motion to satisfy $\boldsymbol v=\boldsymbol 0$ everywhere in the cavity and in the same process, the initial angular momentum of the liquid is transferred to the solid. After this transition, the previous considerations relative to the sensitivity on the initial angular velocity apply.

\subsection{The ``flip-over" effect}
\begin{figure}
\begin{center}
$\nu=0.0375,\, \angvel_0 = [6.2697,0.4109,0],\, E(0)=0$\\
\includegraphics[width=\textwidth]{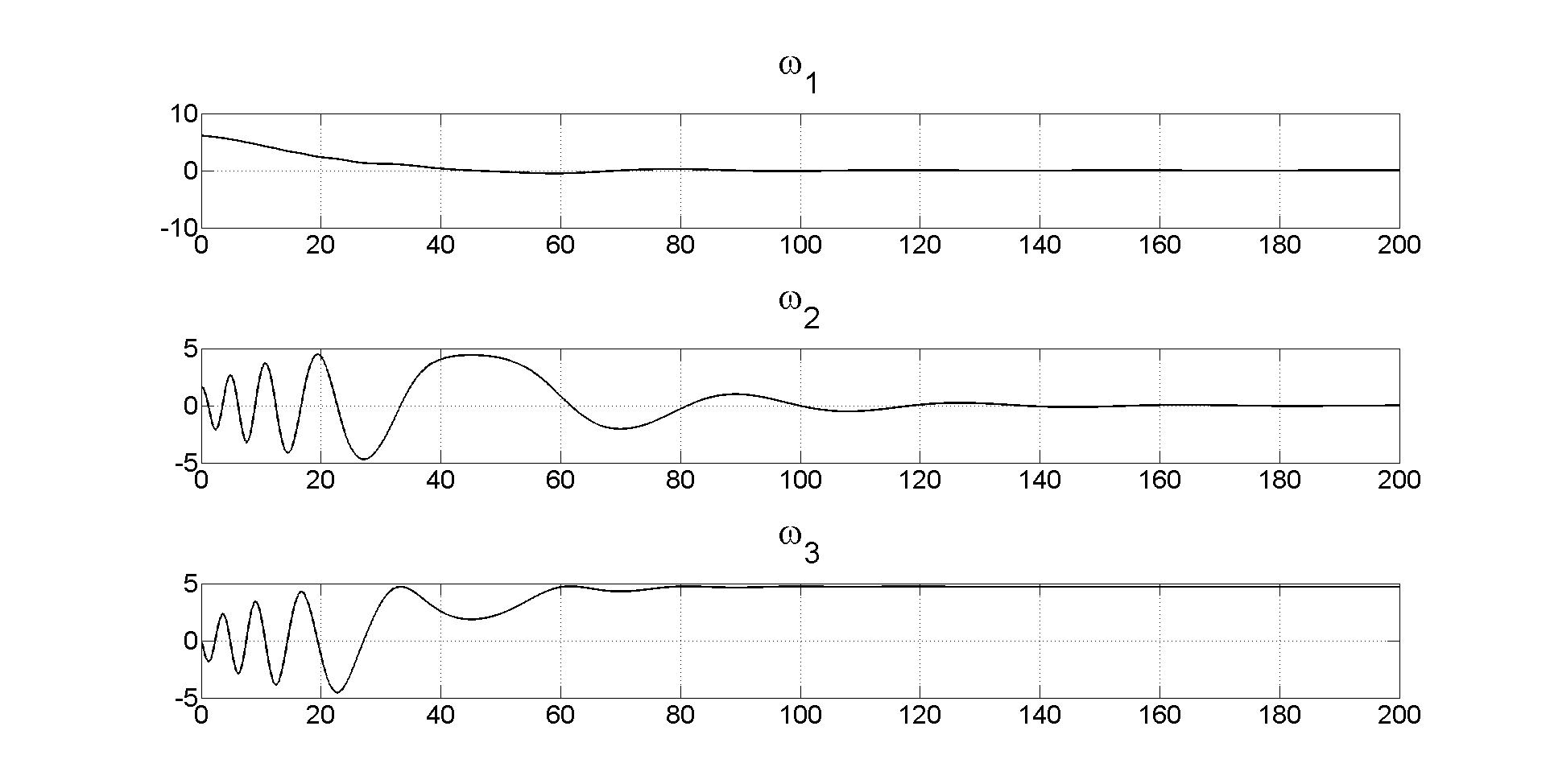}\\

$\nu=0.035,\, \angvel_0 = [6.2697,0.4109,0],\, E(0)=0$\\
\includegraphics[width=\textwidth]{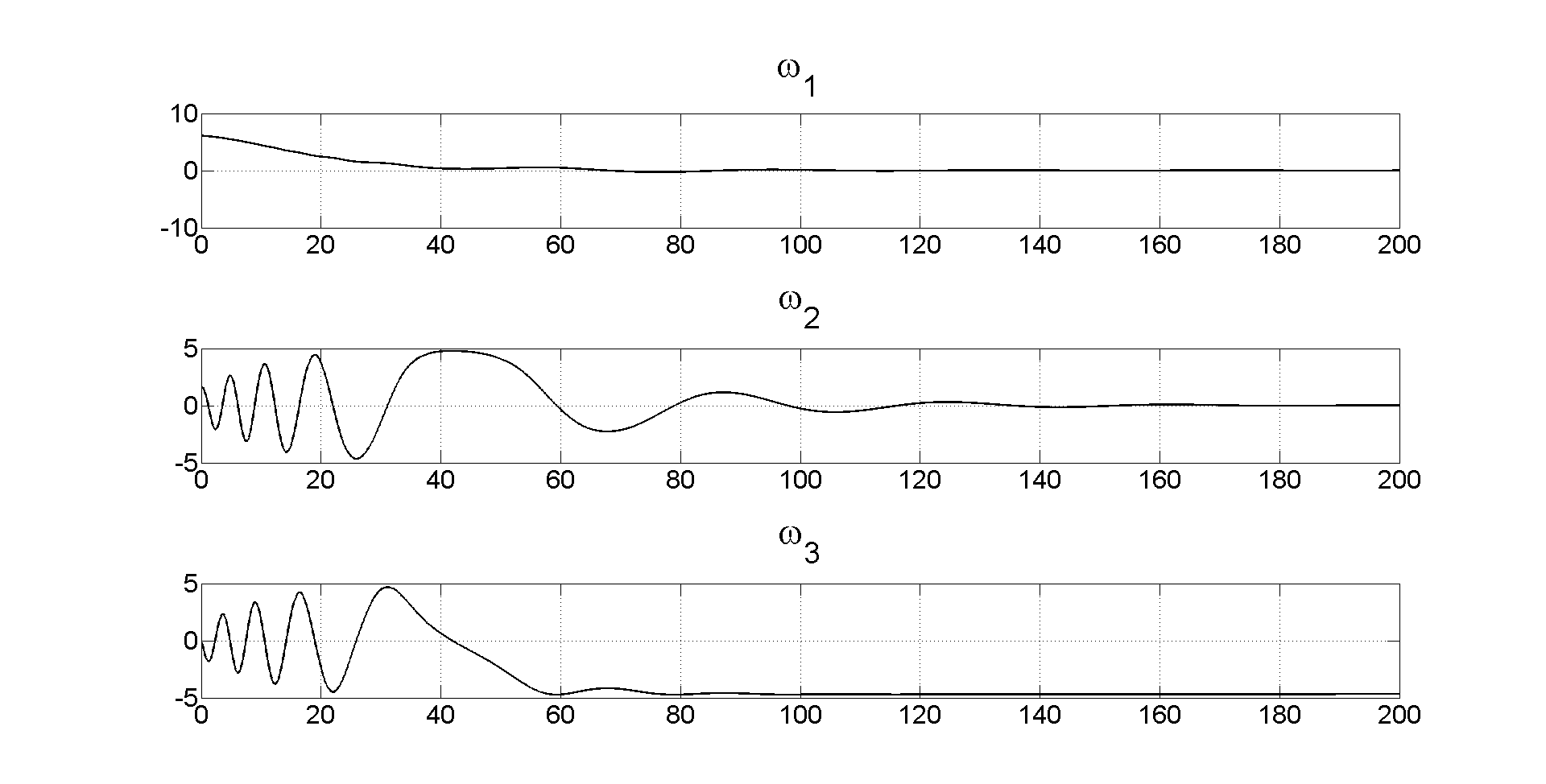}\\

$\nu=0.0325,\, \angvel_0 = [6.2697,0.4109,0],\, E(0)=0$\\
\includegraphics[width=\textwidth]{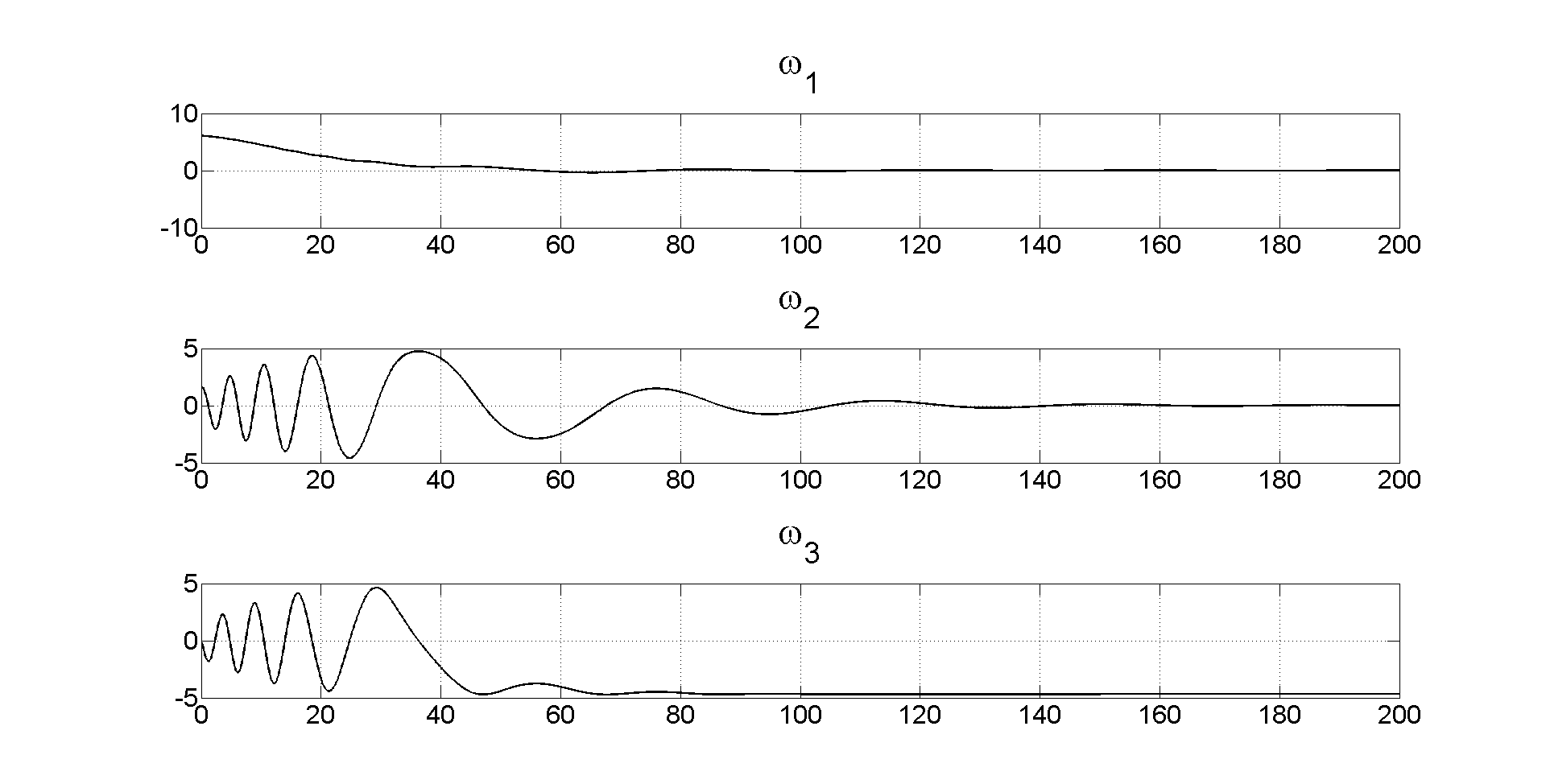}\\
\caption{Visualization of the ``flip-over" effect. The orientation of the final rotation changes, when moving form $\nu=0.0375$ to $\nu=0.035$.}
\label{fig:test4}
\end{center}
\end{figure}

In this section we discuss how, besides the magnitude of the initial rotation along the unstable axes, the viscosity can also trigger the ``flip-over" effect, in some particular configurations. As discussed in Remark 9, if $A \leq B < C$ and $r(0)=0$, the analysis is not sufficient to determine the orientation of the final rotation along the stable axis, namely $\boldsymbol e_3$. The numerical simulations reflect this type of uncertainty and suggest that the determining factor is the liquid viscosity. More precisely, a careful analysis of Table \ref{tab:1} shows that, when $r(0)=0$, changing the viscosity of the liquid not only affects the time to reach equilibrium, but also the orientation of the final rotation, namely $\mathrm{sign}(\overline{r}$). Indeed, jumping from $\nu=0.05$ to $\nu=0.02$, the component $r(t_c)$ changes its sign, while the modulus is almost invariant. The results of Figure \ref{fig:test4} illustrate this effect with more details. For these tests, we consider the case $A=5.54,\, B=6.73,\, C=6.76$ and the initial rotation $\angvel_0 = [6.2697,0.4109,0]$. Figure \ref{fig:test4} shows that, for fixed initial conditions, the sign of $r$ is sensitive to the liquid viscosity. On the basis of numerical experiments, we believe there exists a precise transition point at which the orientation of the rotation is flipped. For the particular configuration considered here, using Figure \ref{fig:test4} we estimate that the transition point is $\nu^* \in (0.0325,0.0375)$.


\renewcommand{\theequation}{A.\arabic{equation}}\setcounter{equation}{0}
\section*{Appendix}\noindent
{\bf Proof of \lemmref{1}.}  Let
\be
r(t)=y(s){\rm e}^{-k(t-s)}+\int_s^t{\rm e}^{-k(t-\tau)}F(\tau)\,d\tau\,.
\eeq{2.2}
By a direct computation one shows that
\be
r(t)=y(s)-k\int_s^tr(\tau)\,d\tau+\int_s^tF(\tau)\,d\tau\,.
\eeq{2.3}
Using the H\"older inequality in the last term of \Eqref{2.2} and recalling that by assumption $y(t)\le M$, for some $M$ independent of $t\ge 0$, and a.a. $t\ge 0$,  we deduce
$$\ba{rl}\smallskip
r(t)\le &M\,{\rm e}^{-k(t-s)}+\left(\Int st{\rm e}^{-kq'(t-\tau)}d\tau\right)^{1/q'}\left(\Int stF^q(\tau)d\tau\right)^{1/q}\\
 \le &M\,{\rm e}^{-k(t-s)}+c\left(\Int st F^q(\tau)d\tau\right)^{1/q}\,.
\ea
$$
In view of the assumption on $F$, from this relation, and from \Eqref{2.2} we show
\be
\lim_{t\to\infty}r(t)=0\ \  \mbox{if $F\not\equiv 0$}\,;\ \ r(t)=y(s)\,{\rm e}^{-k(t-s)}
\ \ \mbox{all $t\ge s$, if $F\equiv 0$}\,.
\eeq{2.4}
We next subtract  \Eqref{2.3} from \Eqref{2.1} side by side, to deduce
\be
y(t)-r(t)\le-k\int_s^t[y(\tau)-r(\tau)]\,d\tau\,.
\eeq{2.5}
Setting
$$
G(t):=\int_s^t[y(\tau)-r(\tau)]\,d\tau\,,
$$
from \Eqref{2.5} we infer
\be
G'(t)\le -k\,G(t)\,, \ \mbox{a.a. $t\ge s$\,,}\ \ G(s)=0\,,
\eeq{2.6}
namely,
 after a simple manipulation, 
$$
\ode{}{t}({\rm e}^{kt}G(t))\le 0\,,
$$ 
which in view of the second relation in \Eqref{2.6} allows us to derive
$$
\int_s^t[y(\tau)-r(\tau)]\,d\tau\le 0\,,\ \ \mbox{for a.a. $s\ge 0$ and all $t\ge s$}\,. 
$$
Therefore,
$$
y(t)\le r(t)\,,\ \ \mbox{for a.a. $t\ge 0$}
$$
and the claimed result follows from the latter relation and \Eqref{2.4}.
\QED
\medskip\par\noindent
{\bf Proof of \lemmref{2}.}
Setting $Y:=y^2$, from (i) we get
\be
Y^\prime\le -2a\,Y+2b\,Y^\beta+F\,,\ \ t\in [t_0,t_1)\,, 
\eeq{2.9}
where $\beta:=(\alpha+1)/2$, $F:=2c\,y$. In view of the second condition in (ii), contradicting \Eqref{2.7} means that there exists $t^*\in (t_0,t_1)$ such that
\be
Y(t)<\delta^2\,,\ \ \mbox{for all $t\in [t_0,t^*)$}\,;\ \ Y(t^*)=\delta^2\,.
\eeq{2.10}
Using this information back in \Eqref{2.9} we find for all $t\in[t_0,t^*)$
$$
Y'(t)\le 2(-a+b\delta^{\alpha-1})Y(t)+F(t)\,,
$$
which in view of the assumptions,  after integration from $t_0$ to $t^*$, furnishes
$$
Y(t^*)< \frac{\delta^2}{2}+\int_{t_0}^{t_1}F(t)\,dt<\delta^2\,. 
$$
However, 
the latter is at odds with \Eqref{2.10}, and we thus conclude the proof of the first part of the 
lemma. In order to show the second part, we observe that from \Eqref{2.7} and \Eqref{2.9} we deduce
$$
Y'\le -2a\,Y+2(b\delta^\alpha+c)y\,,
$$
so that \Eqref{2.8} follows from the first condition in (ii) and \lemmref{1}.
\QED\medskip\par\noindent
{\bf Proof that formulations \Eqref{3.1} and \Eqref{3.5} are equivalent.} We shall show that every (sufficiently smooth) solution to either one satisfies the other as well. We begin to observe that, recalling that $\bfu=\bfv+\bfomega\times\bfx$, equations \Eqref{3.1}$_{1,2}$ and \Eqref{3.5}$_{1,2}$ are clearly equivalent. We next observe that, for any symmetric tensor field $\bfR$ the following identity holds
$$
\int_{\mathscr C}\bfx\times\Div\bfR=\int_{\partial \mathscr C}\bfx\times\bfR\cdot\bfn\,,
$$
where $\bfn$ is the outer unit normal to $\partial\mathscr C$. Therefore, cross-multiplying both sides of \Eqref{3.5}$_1$ by $\bfx$, integrating over $\mathscr C$, and observing that, by \Eqref{3.2} and \Eqref{3.5}$_2$  
$$
\mu\Delta\bfv-\nabla {\rm p}-\rho\bfv\cdot\nabla\bfv=\Div\left[\bfT(\bfv,{\rm p})-\rho\,\bfv\otimes\bfv\right]\,,
$$
we deduce
$$\ba{rl}\smallskip
\ode{}{t}\Int{\mathscr C}{}\rho\bfx\times\bfv +&\Int{\mathscr C}{}\rho\left[2\bfx\times(\bfomega\times\bfv)+\bfx\times(\dot{\bfomega}\times\bfx)\right]\\
&=\Int{\partial\mathscr C}{}\left[\bfx\times\bfT(\bfv,{\sf p})\cdot\bfn-\half\rho\bfx\times\bfn(\bfomega\times\bfx)^2\right]\,.
\ea
$$
Using \Eqref{3.1}$_3$ in the latter displayed equation produces
\be
\ba{rl}\smallskip
\ode{}{t}\Int{\mathscr C}{}\rho\bfx\times\bfv +&\Int{\mathscr C}{}\rho\left[2\bfx\times(\bfomega\times\bfv)+\bfx\times(\dot{\bfomega}\times\bfx)\right]\\
&=-\bfI_{\mathscr B}\cdot\ode{\bfomega}{t}-\bfomega\times(\bfI_{\mathscr B}\cdot\bfomega)-\half\Int{\partial\mathscr C}{}\rho\bfx\times\bfn(\bfomega\times\bfx)^2\,.
\ea
\eeq{A.1}
We now observe the following facts.  
In the first place, using the well-known identities
$$
\int_{\mathscr C}\nabla\times\bfF=\Int{\partial\mathscr C}{}\bfn\times\bfF\,,\ \ 
\ \nabla\times(\varphi\,\bfF)=\varphi\,(\nabla\times\bfF)-\bfF\times\nabla\varphi\,,
$$
we deduce
\be
-\half\Int{\partial\mathscr C}{}\rho\,\bfx\times\bfn(\bfomega\times\bfx)^2=\half\Int{\partial\mathscr C}{}\rho\,\nabla\times[\bfx\,(\bfomega\times\bfx)^2]=-\Int{\mathscr C}{}\rho\,\bfx\times[\bfomega\times(\bfomega\times\bfx)]\,.
\eeq{A.2}
Furthermore, denoting by $\bfI_{\mathscr L}$ the inertia tensor,  with respect to $G$,  of the cavity filled with the liquid, for any $\bfa\in\real^3$ we have
\be
\bfI_{\mathscr L}\cdot\bfa =\Int{\mathscr C}{}\rho\bfx\times({\bfa}\times\bfx)
\,,\ \ \bfa\times(\bfI_{\mathscr L}\cdot\bfa) =\Int{\mathscr C}{}\rho\bfx\times[\bfa\times({\bfa}\times\bfx)]
\,. 
\eeq{A.3}
As a consequence, from \Eqref{A.1}, \Eqref{A.2} and \Eqref{A.3} we conclude
\be
\ode{}{t}\left(\int_{\mathscr C}\rho\,\bfx\times\bfv+\bfI\cdot\bfomega\right)=-\bfomega\times(\bfI\cdot\bfomega)-2\Int{\mathscr C}{}\rho\,\bfx\times(\bfomega\times\bfv)\,,
\eeq{A.4}
where, we recall, $\bfI=\bfI_{\mathscr L}+\bfI_{\mathscr B}$. We next show that
\be
2\int_{\mathscr C}\bfx\times (\bfomega\times\bfv)=\int_{\mathscr C}\bfomega\times (\bfx\times\bfv)\,.
\eeq{A.5}
Let us denote by $\bfL$ the left-hand side of this identity. Since $\Div\bfv=0$, we get
$$
\bfomega\times\bfv=\bfv\cdot\nabla(\bfomega\times\bfx)=\Div [\bfv\otimes(\bfomega\times\bfx)].
$$
Therefore, integrating by parts and taking into account that $\bfv$ vanishes on $\partial\mathscr C$
$$
L_i=2\int_{\mathscr C}\big(\bfx\times\Div [\bfv\otimes(\bfomega\times\bfx)]\big)_i=-2\int_{\mathscr C}\epsilon_{ilk}\epsilon_{kmn}v_l\omega_mx_n\,.
$$
Recalling that $\epsilon_{ilk}\epsilon_{kmn}=\delta_{im}\delta_{ln}-\delta_{in}\delta_{lm}$, we thus deduce
$$
\bfL=-2\int_{\mathscr C}\bfomega(\bfv\cdot\bfx)+2\int_{\mathscr C}(\bfomega\cdot\bfv)\bfx\,,
$$
from which, since $\bfx=\half\nabla(\bfx\cdot\bfx)$, by the assumption on $\bfv$ we conclude
\be
\bfL=2\int_{\mathscr C}(\bfomega\cdot\bfv)\bfx\,.
\eeq{st}
On the other hand, using the identity 
\be
(\bfw\times\bfz)\times\bfs=(\bfw\cdot\bfs)\bfz-(\bfz\cdot\bfs)\bfw
\eeq{st1}
we get
$$
\int_{\mathscr C}\bfx\times(\bfomega\times\bfv)=-\int_{\mathscr C}(\bfx\cdot \bfomega)\bfv+\int_{\mathscr C}
(\bfx\cdot\bfv)\bfomega
=-\int_{\mathscr C}(\bfx\cdot \bfomega)\bfv,
$$
by the assumption on $\bfv$. 
Combining the latter with \Eqref{st} and using  \Eqref{st1} produces
\be
\bfL=\int_{\mathscr C}\bfx(\bfomega\cdot\bfv)-\int_{\mathscr C}(\bfx\cdot \bfomega)\bfv\,
=\int_{\mathscr C}\bfomega \times (\bfx \times \bfv),
\eeq{st2}
the identity \Eqref{A.5} then follows. 

From \Eqref{A.4} and \Eqref{A.5} we obtain \Eqref{3.5}$_3$, and then conclude that every solution to \Eqref{3.1} necessarily satisfies \Eqref{3.5}. Conversely, by retracing backward the calculation just performed, we show that every solution to \Eqref{3.5} also satisfies \Eqref{3.1}, and this concludes the proof of the equivalence of the two formulations. \QED
\medskip\par\noindent
{\em Acknowledgment.} The main results reported in this paper were presented at different scientific venues, including ``International Conference 
 on the Mathematical Fluid Dynamics'' in honor of Professor {\sc Y. Shibata}'s 60th birthday (Waseda, Japan, March 5-9 2013),   
 ``Navier-Stokes Equations in Venice'' (Venice, Italy, April 8-12 2013), and ``Workshop on Navier--Stokes Equations'' (Aachen, Germany, May 21-24 2013). The generous financial support of the organizing Institutions is greatly appreciated.  
This work was also partially supported by the NSF grant DMS-1311983.

Giovanni P. Galdi,\ Giusy Mazzone \& Paolo Zunino
\\
Department of Mechanical Engineering \& Materials Science\\
University of Pittsburgh,\  
Pittsburgh, PA 15261 \\
email: {galdi@pitt.edu, gim20@pitt.edu, paz13@pitt.edu}        
\ed
\begin{thebibliography}{99}
\bibitem{Ball}Ball, J., Continuity Properties and Global Attractors of
Generalized Semiflows and the Navier--Stokes Equations, {\em J. Nonlinear Sci.} {\bf 7}   475--502 (1997)
\bibitem{MR1115205}
Brezzi, F.,  and Fortin M.,
\newblock {\em Mixed and hybrid finite element methods}, volume~15 of {\em
  Springer Series in Computational Mathematics}.
\newblock Springer-Verlag, New York (1991)

\bibitem{Ch}Chernousko, F.L.,  Motion of a Rigid Body with Cavities Containing a Viscous Fluid (1968).
Moscow; NASA Technical Translations, May 1972
\bibitem{DaSl}Dafermos, C.M., and Slemrod, M.,
Asymptotic Behavior of Nonlinear Contraction Semigroups.
{\em J. Functional Analysis}, {\bf 13} 97--106 (1973) 
\bibitem{Gabook}Galdi, G.P., \emph{An Introduction to the Mathematical Theory of the Navier--Stokes Equations:  Steady-State Problems}, Springer-Verlag, New York, Second Edition (2011)
\bibitem{GMZ}Galdi, G.P., Mazzone, G., and Zunino, P., Inertial Motions of a Rigid Body with a Cavity Filled with a Viscous Liquid,
 {\em Comptes Rendus M\'ecanique},  {\bf 341} 760--765 (2013)
\bibitem{MR851383}
Girault, V., and Raviart, P-A.,
\newblock {\em Finite element methods for {N}avier-{S}tokes equations},
  volume~5 of {\em Springer Series in Computational Mathematics}.
\newblock Springer-Verlag, Berlin (1986)
\bibitem{MR1017032}
Gunzburger, M.D.,
\newblock {\em Finite element methods for viscous incompressible flows}.
\newblock Computer Science and Scientific Computing. Academic Press Inc.,
  Boston, MA (1989)
\bibitem{MR2221614}
Hairer, E.,  Lubich, Ch., and  Wanner, G., 
\newblock {\em Geometric numerical integration}, volume~31 of {\em Springer
  Series in Computational Mathematics}.
\newblock Springer-Verlag, Berlin, second edition (2006)
 
\bibitem{Hale}Hale, J.K., {\em Asymptotic Behavior of Dissipative Systems}. Mathematical Surveys and Monographs, Vol. 25. American Mathematical Society, Providence, RI  (1988)
\bibitem{Hough}
Hough, S.S.,  On the Oscillations of a Rotating Ellipsoidal Shell Containing Fluid, {\em Phil. Trans. Roy. Soc. London}, {\bf 186}, 469--506 (1895)
\bibitem{KK}Kopachevsky, N.D., and Krein, S.G., {\em Operator Approach to Linear 
Problems of Hydrodynamics, Vol.2: Nonself-Adjoint Problems for 
Viscous Fluids}, Birkh\"auser Verlag, Basel-Boston-Berlin (2000)
\bibitem{KSY} Kostyuchenko, A.G.,  Shkalikov, A.A., and Yurkin, M.Yu., On the Stability of a Top with a Cavity Filled with a Viscous Fluid, {\em   Funct. Anal. Appl.}, {\bf 32} 100--113 (1998)
\bibitem{Lei}Leimanis, E. {\em The General Problem of the Motion of Coupled Rigid Bodies about a Fixed Point}, Springer Tracts in Natural Philosophy, Vol. 7 (1965)
\bibitem{LK}Leung, A.Y.T., and Kuang, J.L., Chaotic Rotations of a Liquid-Filled Solid,
 {\em J. Sound and Vibr.} {\bf 302} 540--563 (2007)
\bibitem{Mazz}Mazzone, G., A Mathematical Analysis
of the Motion of a Rigid Body
with a Cavity
Containing a Newtonian Fluid, PhD Thesis, Department of Mathematics, Universit\`a del Salento  (2012)
\bibitem{MoRu}Moiseyev, N.N., and Rumyantsev, V.V., {\em Dynamic Stability of Bodies Containing Fluid}, Springer-Verlag, New York (1968)
\bibitem{Prodi}Prodi, G., Teoremi di Tipo Locale per il Sistema di Navier--Stokes e Stabilit\`a delle Soluzioni Stazionarie, {\em
Rend. Sem. Mat. Univ. Padova}, {\bf 32}  374--397 (1962)
\bibitem{MR2454024}
 Roos, H-G.,  Stynes, M., and Tobiska, L.,
\newblock {\em Robust numerical methods for singularly perturbed differential
  equations}, volume~24 of {\em Springer Series in Computational Mathematics}.
\newblock Springer-Verlag, Berlin, second edition  (2008)
\bibitem{Rum}Rumyantsev, V.V., Lyapunov Methods in Analysis of Stability of Motion of Rigid Bodies with Cavities
Filled with Fluid, {\em Izv. Akad. Nauk SSSR, Ser. Mech.}, {\bf 6} 119--140 (1963)
\bibitem{ST} Silvestre, A.L., and Takahashi, T., On the Motion of a Rigid Body with a Cavity Filled with a Viscous Liquid, {\em Proc. Roy. Soc. Edinburgh} Sect. A  {\bf 142}391--423  (2012)
\bibitem{Smirn1}Smirnova, E.P., Stabilization of Free Rotation of an Asymmetric Top with Cavities Entirely Filled with
Fluid, {\em Prikl. Mat. Mekh.}, {\bf 38} 980--985 (1974)
\bibitem{Smirn2}Smirnova, E.P.,  Stability of Free Rotation of a Top Containing a Toroidal Cavity with a Viscous fluid of
Small Viscosity, {\em Mekh. Tverd. Tela}, {\bf  5} 20--25 (1976)
\bibitem{Sob}Sobolev, S.L., Motion of a Symmetric Top with a Cavity Filled with a Fluid, {\em Zh. Prikl. Mekh. Tekhn.
Fiz.}, {\bf 3} 20--55 (1960) 
\bibitem{Stokes0} Stokes, G.,  {\em Mathematical and Physical Papers}, {Vol.  1}, Cambridge (1880) 
\bibitem{Zu}
Zhukovskii, N.Ye., On the Motion of a Rigid Body with Cavities Filled with a 
Homogeneous Liquid Drop, {\em Zh. Fiz.-Khim. Obs.} physics part, {\bf 17} 81--113 (1885);
{\bf 17}  145--199 (1885);
{\bf 17} 231--280 (1885)  
Reprinted in his Selected Works, {\bf 1} (Gostekhizdat, Moscow, 1948) 31--152.




\end{thebibliography}
